\documentstyle[11pt]{article}
\begin{document}
\includegraphics{zibheader.eps}
\vspace*{6.5cm}
\begin{center}
{\Large Wolfram Koepf}\vspace*{2cm}\\
\LARGE{{\bf Algorithms for the}}\\[2mm]
\LARGE{{\bf Indefinite and Definite Summation}}
\end{center}
\vfill
\hrule
\vspace*{3mm}
Preprint SC 94--33 (Dezember 1994)

\thispagestyle{empty}
\setcounter{page}{0}
\eject
\hoffset -1cm
\voffset -2cm
\footskip=1cm
\parindent=0pt

\newcommand{\Maple}{{\sc Maple}}
\newcommand{\Macsyma}{{\sc Macsyma}}
\newcommand{\Mathematica}{{\sc Mathematica}}
\newcommand{\Reduce}{{\sc Reduce}}
\def\lcm{\mathop{\rm lcm}\nolimits}
\newcommand{\PSbild}[4]{\begin{center}\vspace*{#1}
\makebox[#2][l]{\includegraphics{#3}}\par
\begin{figure}[h]
\caption{#4}
\label{#3}
\end{figure}
\end{center}
}

\def\arcsin{\mathop{\rm arcsin}\nolimits}
\def\arccos{\mathop{\rm arccos}\nolimits}
\def\arctan{\mathop{\rm arctan}\nolimits}
\def\arccot{\mathop{\rm arccot}\nolimits}
\def\arcsec{\mathop{\rm arcsec}\nolimits}
\def\arccsc{\mathop{\rm arccsc}\nolimits}
\def\arcsinh{\mathop{\rm arsinh}\nolimits}
\def\arccosh{\mathop{\rm arcosh}\nolimits}
\def\arctanh{\mathop{\rm artanh}\nolimits}
\def\arccoth{\mathop{\rm arcoth}\nolimits}
\def\arcsech{\mathop{\rm arsech}\nolimits}
\def\arccsch{\mathop{\rm arcsch}\nolimits}
\def\sech{\mathop{\rm sech}\nolimits}
\def\csch{\mathop{\rm csch}\nolimits}
\def\erf{\mathop{\rm erf}\nolimits}
\newcommand{\one}[2]{\\[1mm] \noindent \begin{tabular}
{p{.5cm}p{12.5cm}}
\mbox{\rm (#1)}&
{$\displaystyle{#2}$}
\end{tabular}}
\newcommand{\two}[4]{\\[1mm] \noindent \begin{tabular}
{p{.5cm}p{6.5cm}p{.5cm}p{6.5cm}}
\mbox{\rm (#1)}&
{$\displaystyle{#2}$}&
\mbox{\rm (#3)}&
{$\displaystyle{#4}$}
\end{tabular}}
\newcommand{\three}[6]{\\[1mm] \noindent \begin{tabular}
{p{.5cm}p{4cm}p{.5cm}p{4cm}p{.5cm}p{4cm}}
\mbox{\rm (#1)}&
{$\displaystyle{#2}$}&
\mbox{\rm (#3)}&
{$\displaystyle{#4}$}&
\mbox{\rm (#5)}&
{$\displaystyle{#6}$}
\end{tabular}}
\newcommand{\four}[8]{\\[1mm] \noindent \begin{tabular}
{p{.5cm}p{2.75cm}p{.5cm}p{2.75cm}p{.5cm}p{2.75cm}p{.5cm}p{2.75cm}}
\mbox{\rm (#1)}&
{$\displaystyle{#2}$}&
\mbox{\rm (#3)}&
{$\displaystyle{#4}$}&
\mbox{\rm (#5)}&
{$\displaystyle{#6}$}&
\mbox{\rm (#7)}&
{$\displaystyle{#8}$}
\end{tabular}}

\newcommand{\RR}{{\rm I\! R}}
\newcommand{\DD}{{\rm I\! D}} 
\newcommand{\CC}{{\; \rm l\!\!\!        C}}
\newcommand{\NN}{{\rm I\!I\!\! N}}
\newcommand{\ZZ} {{\rm {\mbox{\protect\makebox[.15em][l]{\small \sf Z}\small \sf
 Z}}}}
\newcommand{\C}{{\rm {\mbox{C{\llap{{\vrule height1.52ex}\kern.4em}}}}}}
\newcommand{\N} {{\rm {\mbox{\protect\makebox[.15em][l]{I}N}}}}
\newcommand{\Q} {{\rm {\mbox{Q{\llap{{\vrule height1.5ex}\kern.5em}}}}}}
\newcommand{\R} {{\rm {\mbox{\protect\makebox[.15em][l]{I}R}}}}
\newcommand{\D} {{\rm {\mbox{\protect\makebox[.15em][l]{I}D}}}}
\newcommand{\K} {{\rm {\mbox{\protect\makebox[.15em][l]{I}K}}}}
\newcommand{\Z} {{\rm {\mbox{\protect\makebox[.2em][l]{\sf Z}\sf Z}}}}
\newcommand{\Pt}{\widetilde{P}}
\newcommand{\Cd}{\widehat{{\rm \;l\!\!\! C}}}
\newcommand{\bq}[1]{|#1|^2}
\newcommand{\vf}[1]{(1-\bq{#1})}
\newcommand{\vfq}[1]{\vf{#1}^2}
\newcommand{\Bf}[1]{\vf{#1}\left|\frac{f''}{f'}(#1)\right|}
\newcommand{\Nf}[1]{\vfq{#1}|S_f(#1)|}
\newcommand{\Kf}[1]{\left|-\kon{#1}+\frac{1}{2}\vf{#1}\frac{f''}{f'}(#1)\right|}
\newcommand{\ed}[1]{\frac{1}{#1}}
\newcommand{\aut}[1]{\frac{{\textstyle{z+#1}}}{{\textstyle{1+\kon{#1}z}}}}
\newcommand{\au}[2]{#1\aut{#2}}
\newcommand{\subs}[2]{\left. \makebox{\rule{0in}{2.5ex}} #1 \rb_{#2}}
\newcommand{\subst}[3]{\left. \makebox{\rule{0in}{2.5ex}} #1 \rb_{#2}^{#3}}
\newcommand{\funkdef}[3]{\left\{\!\!\!\begin{array}{cc}
                                #1 & \!\!\!\mbox{\rm{if} $#2$ } \\
                                #3 & \!\!\!\mbox{\rm{otherwise}}
                                \end{array}  
                         \right.}       
\newcommand{\funkdeff}[4]{\left\{\begin{array}{ccc}
                                 #1 && \mbox{\rm{if} $#2$ } \\[3mm]
                                 #3 && \mbox{\rm{if} $#4$ } 
                                 \end{array}
                          \right.}
\newcommand{\funkdefff}[6]{\left\{\begin{array}{ccc}
                                 #1 && \mbox{\rm{if} $#2$ } \\
                                 #3 && \mbox{\rm{if} $#4$ } \\
                                 #5 && \mbox{\rm{if} $#6$ }
                                 \end{array}
                           \right.}
\newcommand{\funkdeffff}[8]{\left\{\begin{array}{ccc}
                                 #1 && \mbox{\rm{if} $#2$ } \\
                                 #3 && \mbox{\rm{if} $#4$ } \\
                                 #5 && \mbox{\rm{if} $#6$ } \\
                                 #7 && \mbox{\rm{if} $#8$ }
                                 \end{array}
                          \right.}
\newcommand{\ueber}[2]{
                       \Big( \!
                       {{\small
                       \begin{array}{c}
                          #1\\
                          #2
                          \end{array}
                       }}
                       \! \Big) }
\newcommand{\function}[4]{
                          \begin{array}{rcl}#1&\pf&#2\\
                          #3&\mapsto &#4
                          \end{array} }
\newcommand{\pr}{\vspace{-2mm}\absatz{{\sl Proof:}}\hspace{5mm}}
\newcommand{\eop}{\hfill$\Box$\par\medskip\noindent}

\newenvironment{proof}{\pr}{\eop}

\newcommand{\absatz}{\par\medskip\noindent}
\renewcommand{\Re}{{\rm Re\:}}
\renewcommand{\Im}{{\rm Im\:}}
\newcommand{\co}{{\rm co\:}}
\newcommand{\coq}{\overline{{\rm co}} \:}
\newcommand{\ex}{{\rm E\:}}
\newcommand{\In}{\in}
\newcommand{\ro}{\varrho}
\newcommand{\om}{\omega}
\newcommand{\al}{\alpha}
\newcommand{\bb}{\beta}
\newcommand{\la}{\lambda}
\newcommand{\eps}{\varepsilon}
\newcommand{\ph}{\varphi}
\renewcommand{\phi}{\varphi}
\newcommand{\si}{\sigma}
\newcommand{\ka}{\varkappa}
\newcommand{\th}{\theta}
\newcommand{\g}{\gamma}
\newcommand{\de}{\partial}
\newcommand{\fD}{f(\D)}
\newcommand{\sumi}{\sum\limits_{k=0}^{\infty}}
\newcommand{\sumei}{\sum\limits_{k=1}^{\infty}}
\newcommand{\union}{\bigcup\limits_{k=1}^{n}}
\newcommand{\sumn}{\sum\limits_{k=1}^{n}}
\newcommand{\prodn}{\prod\limits_{k=1}^{n}}
\newcommand{\intd}{\int\limits_{\de\DD}}
\newcommand{\menge}[3]{\left\{#1 \In #2 \; \lb \; #3 \right. \right\} }
\newcommand{\mk}{\mu_{k}}
\newcommand{\xk}{x_k}
\newcommand{\yk}{y_k}
\newcommand{\xn}{x_n}
\newcommand{\yn}{y_n}
\newcommand{\ak}{\al_k}
\newcommand{\bk}{\bb_k}
\newcommand{\kn}{\mbox{$(k=1, \ldots ,n)$}}
\newcommand{\kno}{\mbox{$k=1, \ldots ,n$}}
\newcommand{\sub}{\prec}
\newcommand{\ld}[1]{\frac{f''}{f'}(#1)} 
\newcommand{\limn}{\lim\limits_{n\rightarrow\infty}}
\newcommand{\lsz}{\limsup\limits_{z\rightarrow\de\DD}}
\newcommand{\limr}{\lim\limits_{r\rightarrow 1}}
\newcommand{\liz}{\liminf\limits_{z\rightarrow\de\DD}}
\newcommand{\supD}[1]{\sup\limits_{#1\In \DD}}
\newcommand{\infD}[1]{\inf\limits_{#1\In \DD}}
\newcommand{\maxn}{\max\limits_{1 \leqq k \leqq n}}
\newcommand{\minn}{\min\limits_{1 \leqq k \leqq n}}
\newcommand{\ord}{{\rm ord\:}}
\newcommand{\gleich}[1]{\stackrel{{{\rm (#1)}}}{\leq}}
\newcommand{\folgt}[1]{\stackrel{{{\rm (#1)}}}{\Pf}}
\newcommand{\gegen}[1]{\stackrel{{{\rm (#1)}}}{\pf}}
\newcommand{\nach}[2]{(#1)$\dpf$(#2):}
\newcommand{\Exp}{\subset\!\subset}
\newcommand{\lleq}{\stackrel{_{{\scriptscriptstyle \Exp}}}
           {_{{\scriptscriptstyle \sim}}}}
\newcommand{\Sub}{{\rm Sub\:}}
\newcommand{\norm}[2]{\frac{#1\circ #2-#1\circ #2(0)}{(#1\circ #2)'(0)}} 
\renewcommand{\dim}{{\rm{dim}}_{_{H^{p}}}}
\renewcommand{\ll}{<\!<}
\newcommand{\1}{{\bf{1}}}
\newcommand{\2}{{\bf{2}}}
\newcommand{\4}{{\bf{4}}}
\newcommand{\5}{{\bf{5}}}
\newcommand{\6}{{\bf{6}}}
\newcommand{\7}{{\bf{7}}}
\newcommand{\8}{{\bf{8}}}
\newcommand{\9}{{\bf{9}}}
\newcommand{\0}{{\bf{0}}}


\newcommand{\luc}{locally uniform convergence}
\newcommand{\hp}{half\-plane}
\newcommand{\sq}{sequence}
\newcommand{\an}{analytic}
\newcommand{\af}{analytic\ function}
\newcommand{\Ne}{Ne\-ha\-ri\ ex\-pression}
\newcommand{\Ke}{Koe\-be\ ex\-pression}
\newcommand{\Be}{Becker\ ex\-pression}
\newcommand{\fc}{function}
\newcommand{\uv}{univalent}
\newcommand{\uf}{univalent\ function}
\newcommand{\SC}{Schwarz-Chri\-stof\-fel}
\newcommand{\pol}{polyg\-on}
\newcommand{\cv}{con\-vex}
\newcommand{\ctc}{close-to-con\-vex}
\newcommand{\Ck}{Ca\-ra\-th\'eo\-dory\ kernel}
\newcommand{\st}{sector}
\newcommand{\sm}{similar}
\newcommand{\br}{boundary\ rotation}
\newcommand{\bbr}{bounded\ \br}
\newcommand{\KM}{Krein-Mil\-man}


\newcommand{\til}{\widetilde}
\newcommand{\pf}{\rightarrow}
\newcommand{\Pf}{\;\;\;\longrightarrow\;\;\;}
\newcommand{\dpf}{\Rightarrow}
\newcommand{\Dpf}{\;\;\;\Longrightarrow\;\;\;}
\newcommand{\kon}{\overline}
\newcommand{\be}{\begin{equation}}
\newcommand{\ee}{\end{equation}}
\newcommand{\bea}{\begin{eqnarray}}
\newcommand{\eea}{\end{eqnarray}}
\newcommand{\beao}{\begin{eqnarray*}}
\newcommand{\eeao}{\end{eqnarray*}}
\newcommand{\lequiv}{\;\;\;\Longleftrightarrow\;\;\;}
\newcommand{\gl}{\;\leftrightarrow\;}
\newcommand{\lk}{\left(}
\newcommand{\rk}{\right)}
\newcommand{\lb}{\left|}
\newcommand{\rb}{\right|}


\newcommand{\bT}{\begin{theorem}}
\newcommand{\eT}{\end{theorem}}
\newcommand{\bL}{\begin{lemma}}
\newcommand{\eL}{\end{lemma}}
\newcommand{\bC}{\begin{corollary}}  
\newcommand{\eC}{\end{corollary}}
\newcommand{\bt}{\begin{tabbing} 12345 \= \kill}
\newcommand{\et}{\end{tabbing}}


\hyphenation{qua-si-disk qua-si-circle 
non-smooth
pa-ram-e-trized pa-ram-e-tri-zation
geo-met-ric
}


\newcommand{\bbegin}{{\bf{begin}}}
\newcommand{\eend}{{\bf{end}}}
\newcommand{\iif}{{\bf{if}}}
\newcommand{\tthen}{{\bf{then}}}
\newcommand{\wwhile}{{\bf{while}}}
\newcommand{\ddo}{{\bf{do}}}
\newcommand{\ffor}{{\bf{for}}}
\newcommand{\sstep}{{\bf{step}}}
\newcommand{\llet}{{\bf{let}}}
\newcommand{\pprocedure}{{\bf{procedure}}}
\newcommand{\aand}{{\bf{and}}}
\newcommand{\nnot}{{\bf{not}}}
\newcommand{\oor}{{\bf{or}}}
\newcommand{\lllet}{{\bf{let}}}
\newcommand{\eexit}{{\bf{exit}}}
\newcommand{\rreturn}{{\bf{return}}}
\newcommand{\uuntil}{{\bf{until}}}

\newcommand{\abs}{\\[3mm]}
\newcommand{\leqq}{\leq}
\newcommand{\geqq}{\geq}

\newtheorem{theorem}{Theorem}
\newtheorem{algorithm}{Algorithm}
\newtheorem{lemma}{Lemma}
\newtheorem{corollary}{Corollary}
\newtheorem{remark}{Remark}
\newtheorem{example}{Example}

\newcommand{\ded}[1]{\frac{\de}{\de #1}}
\newcommand{\dedn}[2]{\frac{\de^{#2}}{\de #1^{#2}}}
\def\airyai{\mathop{\rm Ai}\nolimits}
\def\airybi{\mathop{\rm Bi}\nolimits}
\def\erfc{\mathop{\rm erfc}\nolimits}
\def\erf{\mathop{\rm erf}\nolimits}

\begin{center}
\LARGE{{\bf Algorithms for the}}\\
\LARGE{{\bf Indefinite and Definite Summation}}\\[3mm]
{\Large Wolfram Koepf}\\
{\large \tt koepf@zib-berlin.de}
\vspace*{2cm}
\end{center}


\label{firstpage}

\begin{abstract}
The celebrated Zeilberger algorithm which finds holonomic recurrence
equations for definite sums of hypergeometric terms $F(n,k)$
is extended to certain nonhypergeometric terms. An expression 
$F(n,k)$ is called hypergeometric term if both $F(n+1,k)/F(n,k)$ and
$F(n,k+1)/F(n,k)$ are rational functions. Typical examples
are ratios of products of exponentials, factorials, $\Gamma$ function terms, binomial
coefficients, and Pochhammer symbols that are integer-linear with
respect to $n$ and $k$ in their arguments.

We consider the more general case of ratios of products of exponentials, 
factorials, $\Gamma$ 
function terms, binomial coefficients, and Pochhammer symbols that are 
rational-linear with
respect to $n$ and $k$ in their arguments, and present an extended version
of Zeilberger's algorithm for this case, using an extended version
of Gosper's algorithm for indefinite summation.

In a similar way the Wilf-Zeilberger method of rational function
certification of integer-linear hypergeometric identities is extended
to rational-linear hypergeometric identities.

The given algorithms on definite summation
apply to many cases in the literature to which neither
the Zeilberger approach nor the Wilf-Zeilberger method is applicable. 
Examples of this type are given by theorems of Watson and Whipple, 
and a large list of identities 
(``Strange evaluations of hypergeometric series'')
that were studied by Gessel and Stanton. It turns out that with our
extended algorithms practically all hypergeometric identities in the
literature can be verified.

Finally we show how the algorithms can be used to generate new identities.

\Reduce\ and \Maple\ implementations of the given algorithms can be obtained 
from the author, many results of which are presented in the paper.
\end{abstract}

\section{Hypergeometric identities}
\noindent
In this paper we deal with hypergeometric identities. As usual, the
notation of the generalized hypergeometric function $_{p}F_{q}$ defined by
\be
_{p}F_{q}\left.\lk\begin{array}{cccc}
a_{1}&a_{2}&\cdots&a_{p}\\
b_{1}&b_{2}&\cdots&b_{q}\\
            \end{array}\rb x\rk
:=
\sumi A_k\,x^{k}=
\sumi \frac
{(a_{1})_{k}\cdot(a_{2})_{k}\cdots(a_{p})_{k}}
{(b_{1})_{k}\cdot(b_{2})_{k}\cdots(b_{q})_{k}\,k!}x^{k}
\label{eq:coefficientformula}
\ee
is used, $(a)_{k}=\frac{\Gamma(a+k)}{\Gamma(k)}$ denoting the
{\sl Pochhammer symbol} or {\sl shifted factorial.} 
The numbers $a_{k}$ are called the {\sl upper}, and $b_{k}$ the {\sl lower
parameters} of $_{p} F_{q}$. 

The coefficients $A_{k}$ of the generalized hypergeometric function have the
rational term ratio 
\[
\frac{A_{k+1}}{A_{k}}=
\frac{(k+a_1)\cdot(k+a_2)\cdots (k+a_p)}{(k+b_1)\cdot(k+b_2)\cdots
(k+b_q)(k+1)}\quad(k\in\N)
\;.
\]
If an expression $A_{k}$ has a rational term ratio
$A_{k+1}/{A_k}$, we call $A_{k}$ a {\sl hypergeometric term}
or {\sl closed form}.
Note that any hypergeometric term essentially has a representation 
as the ratio of shifted factorials (over $\C$), and its generating function 
is connected with a generalized hypergeometric series.

{\small
\begin{table}[hbtp]
\caption{Bailey's hypergeometric database}
\label{table:Bailey}
\vspace*{3mm}
\begin{tabular}{rlp{11cm}}
page $\!\!\!$& Theorem & Identity\\[2mm] \hline \\[-2mm]
2--3 & 
\multicolumn{1}{l}
{\hspace*{-2mm}\begin{tabular}{l}Vandermonde\\ Gau{\ss}\end{tabular}}
&
$\displaystyle
_2 F_1\left.
\!\!
\left(
\!\!\!\!
\begin{array}{c}
\multicolumn{1}{c}{\begin{array}{c}
a\;, b
\end{array}}\\[1mm]
\multicolumn{1}{c}{\begin{array}{c}
c
            \end{array}}\end{array}
\!\!\!\!
\right| 1\right)
=\frac{(c-b)_{-a}}
{(c)_{-a}}=
\frac{\Gamma(c)\Gamma(c-a-b)}{\Gamma(c-a)\Gamma(c-b)}
$
\\[5mm]
9 & Saalsch\"utz & 
$\displaystyle
_3 F_2\left.
\!\!
\left(
\!\!\!\!
\begin{array}{c}
\multicolumn{1}{c}{\begin{array}{c}
a\;, b\;, -n
\end{array}}\\[1mm]
\multicolumn{1}{c}{\begin{array}{c}
c\;,1+a+b-c-n
            \end{array}}\end{array}
\!\!\!\!
\right| 1\right)
=\frac{(c-a)_n\,(c-b)_{n}}
{(c)_n\,(c-a-b)_{n}}
$
\\[5mm]
9 & Kummer & 
$\displaystyle
_2 F_1\left.
\!\!
\left(
\!\!\!\!
\begin{array}{c}
\multicolumn{1}{c}{\begin{array}{c}
a\;, b
\end{array}}\\[1mm]
\multicolumn{1}{c}{\begin{array}{c}
1+a-b
            \end{array}}\end{array}
\!\!\!\!
\right| -1\right)
=\frac{(1+a)_{-b}}
{(1+a/2)_{-b}}=
\frac{\Gamma(1+a-b)\Gamma(1+a/2)}{\Gamma(1+a)\Gamma(1+a/2-b)}
$
\\[5mm]
11 & Gau{\ss} & 
$\displaystyle
_2 F_1\left.
\!\!
\left(
\!\!\!\!
\begin{array}{c}
\multicolumn{1}{c}{\begin{array}{c}
a\;, b
\end{array}}\\[1mm]
\multicolumn{1}{c}{\begin{array}{c}
(a+b+1)/2
            \end{array}}\end{array}
\!\!\!\!
\right| \frac{1}{2}\right)
=\frac{\Gamma(1/2)\Gamma((a+b+1)/2)}
{\Gamma((a+1)/2)\Gamma((b+1)/2)}
$
\\[5mm]
11 & Bailey & 
$\displaystyle
_2 F_1\left.
\!\!
\left(
\!\!\!\!
\begin{array}{c}
\multicolumn{1}{c}{\begin{array}{c}
a\;, 1-a
\end{array}}\\[1mm]
\multicolumn{1}{c}{\begin{array}{c}
c
            \end{array}}\end{array}
\!\!\!\!
\right| \frac{1}{2}\right)
=\frac{\Gamma(c/2)\Gamma((c+1)/2)}
{\Gamma((a+c)/2)\Gamma((1-a+c)/2)}
$
\\[5mm]
13 & Dixon & 
$\displaystyle
_3 F_2\left.
\!\!
\left(
\!\!\!\!
\begin{array}{c}
\multicolumn{1}{c}{\begin{array}{c}
a\;, b\;, c
\end{array}}\\[1mm]
\multicolumn{1}{c}{\begin{array}{c}
1+a-b\;,1+a-c
            \end{array}}\end{array}
\!\!\!\!
\right| 1\right)
=\frac{(1+a)_{-c}\,(1+a/2-b)_{-c}}
{(1+a/2)_{-c}\,(1+a-b)_{-c}}
=$\newline
\hspace*{1.5cm}$\displaystyle
\frac{\Gamma(1+a/2)\Gamma(1+a-b)\Gamma(1+a-c)\Gamma(1+a/2-b-c)}
{\Gamma(1+a)\Gamma(1+a/2-b)\Gamma(1+a/2-c)\Gamma(1+a-b-c)}
$
\\[5mm]
16 & 
\multicolumn{1}{l}
{\hspace*{-2mm}\begin{tabular}{l}Watson\\ Whipple\end{tabular}}
&
$\displaystyle
_3 F_2\left.
\!\!
\left(
\!\!\!\!
\begin{array}{c}
\multicolumn{1}{c}{\begin{array}{c}
a\;, b\;, c
\end{array}}\\[1mm]
\multicolumn{1}{c}{\begin{array}{c}
(a+b+1)/2\;,2c
            \end{array}}\end{array}
\!\!\!\!
\right| 1\right)
=
\frac{\Gamma(\frac{1}{2})\Gamma(\frac{1+2c}{2})
\Gamma(\frac{1+a+b}{2})
\Gamma(\frac{1-a-b+2c}{2})}
{\Gamma(\frac{1+a}{2})
\Gamma(\frac{1+b}{2})\Gamma(\frac{1-a+2c}{2})
\Gamma(\frac{1-b+2c}{2})}
$
\\[5mm]
16 & Whipple &
$\displaystyle
_3 F_2\left.
\!\!
\left(
\!\!\!\!\!\!
\begin{array}{c}
\multicolumn{1}{c}{\begin{array}{c}
a, 1-a, c
\end{array}}\\[1mm]
\multicolumn{1}{c}{\begin{array}{c}
e,1+2c-e
            \end{array}}\end{array}
\!\!\!\!\!\!
\right| 1\right)
=
\frac{\pi2^{1-2c}\Gamma(e)\Gamma(1+2c-e)}
{\Gamma(\frac{a+e}{2})
\Gamma(\frac{a+1+2c-e}{2})\Gamma(\frac{1-a+e}{2})
\Gamma(\frac{2+2c-a-e}{2})}
$
\\[5mm]
26 & 
\multicolumn{1}{l}
{\hspace*{-2mm}\begin{tabular}{l}Dougall's\\ Theorem\end{tabular}}
&
$\displaystyle
_7 F_6\left.
\!\!
\left(
\!\!\!\!
\begin{array}{c}
\multicolumn{1}{c}{\begin{array}{c}
a\;, 1\!+\!a/2\;, b\;, c\;, d\;, 1+2a-b-c-d+n\;, -n
\end{array}}\\[1mm]
\multicolumn{1}{c}{\begin{array}{c}
a/2\;,1\!+\!a\!-\!b\;,1\!+\!a\!-\!c\;, 1\!+\!a\!-\!d\;,
b\!+\!c\!+\!d\!-\!a\!-\!n\;,
1\!+\!a\!+\!n
            \end{array}}\end{array}
\!\!\!\!
\right| 1\right)
=$
\newline
\hspace*{1.5cm}$\displaystyle
\frac{(1+a)_n\,(1+a-b-c)_{n}\,(1+a-b-d)_n\,(1+a-c-d)_n}
{(1+a-b)_n(1+a-c)_n\,(1+a-d)_{n}\,(1+a-b-c-d)_n}
$
\\[5mm]
25/27 & Dougall &
$\displaystyle
_5 F_4\left.
\!\!
\left(
\!\!\!\!
\begin{array}{c}
\multicolumn{1}{c}{\begin{array}{c}
a, 1\!+\!a/2, c, d, e
\end{array}}\\[1mm]
\multicolumn{1}{c}{\begin{array}{c}
a/2,1\!+\!a\!-\!c, 1\!+\!a\!-\!d, 1\!+\!a\!-\!e
            \end{array}}\end{array}
\!\!\!\!
\right| 1\right)
=\frac{(1\!+\!a)_{-e} (1\!+\!a\!-\!c\!-\!d)_{-e}}
{(1\!+\!a\!-\!c)_{-e} (1\!+\!a\!-\!d)_{-e}}
=$\newline
\hspace*{1.5cm}$\displaystyle
\frac{\Gamma(1+a-c)\Gamma(1+a-d)\Gamma(1+a-e)\Gamma(1+a-c-d-e)}
{\Gamma(1+a)\Gamma(1+a-d-e)\Gamma(1+a-c-e)\Gamma(1+a-c-d)}
$
\\[5mm]
28 & Whipple &
\mbox{
$\displaystyle
_4 F_3\left.
\!\!
\left(
\!\!\!\!\!\!
\begin{array}{c}
\multicolumn{1}{c}{\begin{array}{c}
a, 1\!+\!a/2, d, e
\end{array}}\\[1mm]
\multicolumn{1}{c}{\begin{array}{c}
a/2, 1\!+\!a\!-\!d, 1\!+\!a\!-\!e
            \end{array}}\end{array}
\!\!\!\!\!
\right|\! -1\!\right)
\!=\!\frac{(1\!+\!a)_{-e}}{(1\!+\!a\!-\!d)_{-e}}
\!=\!\frac{\Gamma(1\!+\!a\!-\!d)\Gamma(1\!+\!a\!-\!e)}
{\Gamma(1\!+\!a)\Gamma(1\!+\!a\!-\!d\!-\!e)}
$
}
\\[5mm]
30 & Bailey &
$\displaystyle
_3 F_2\left.
\!\!
\left(
\!\!\!\!
\begin{array}{c}
\multicolumn{1}{c}{\begin{array}{c}
a\;, 1+a/2\;, -n
\end{array}}\\[1mm]
\multicolumn{1}{c}{\begin{array}{c}
a/2\;, w
            \end{array}}\end{array}
\!\!\!\!
\right| 1\right)
=\frac{(w-a-1-n)(w-a)_{n-1}}
{(w)_n}
$
\\[5mm]
30 & Bailey &
$\displaystyle
_3 F_2\left.
\!\!
\left(
\!\!\!\!
\begin{array}{c}
\multicolumn{1}{c}{\begin{array}{c}
a\;, b\;, -n
\end{array}}\\[1mm]
\multicolumn{1}{c}{\begin{array}{c}
1+a-b\;, 1+2b-n
            \end{array}}\end{array}
\!\!\!\!
\right| 1\right)
=\frac{(a-2b)_n(1+a/2-b)_n(-b)_n}
{(1+a-b)_n(a/2-b)_n(-2b)_n}
$
\\[5mm]
30 & Bailey &
$\displaystyle
_4 F_3\left.
\!\!
\left(
\!\!\!\!
\begin{array}{c}
\multicolumn{1}{c}{\begin{array}{c}
a\;, 1+a/2\;, b\;, -n
\end{array}}\\[1mm]
\multicolumn{1}{c}{\begin{array}{c}
a/2\;, 1+a-b\;, 1+2b-n
            \end{array}}\end{array}
\!\!\!\!
\right| 1\right)
=\frac{(a-2b)_n(-b)_n}
{(1+a-b)_n(-2b)_n}
$
\\[5mm]
30 & Bailey &
$\displaystyle
_4 F_3\left.
\!\!
\left(
\!\!\!\!\!\!\!
\begin{array}{c}
\multicolumn{1}{c}{\begin{array}{c}
a, 1\!+\!a/2, b , -n
\end{array}}\\[1mm]
\multicolumn{1}{c}{\begin{array}{c}
a/2, 1\!+\!a\!-\!b, 2\!+\!2b\!-\!n
            \end{array}}\end{array}
\!\!\!\!\!\!\!
\right| 1\right)
=\frac{(a\!-\!2b\!-\!1)_n(1/2\!+\!a/2\!-\!b)_n(-b\!-\!1)_n}
{(1\!+\!a\!-\!b)_n(a/2\!-\!b\!-\!1/2)_n(-2b\!-\!1)_n}
$
\end{tabular}
\end{table}
}\noindent
The classical reference concerning generalized hypergeometric series
is the book of Bailey (1935) containing a huge amount of
relations between hypergeometric series some of which represent
the value of certain hypergeometric functions at a special point
(mostly $x=1$ or $x=-1$) by a single hypergeometric term. We will be
concerned with this type of identities, and Table~\ref{table:Bailey}
is a complete list of all such hypergeometric identities found in Bailey's 
book.

Here $n\in\N$ is assumed to represent a positive integer so that the
hypergeometric series with upper parameter $-n$ are terminating.
All other parameters involved represent arbitrary complex variables
such that none of the lower parameters corresponds to a negative integer.

With a method due to Wilf and Zeilberger, and with an algorithm of Zeilberger,
many of these hypergeometric identities can be checked. It turns out, however,
that for some of these identities both methods fail. We give extensions
of both the Wilf-Zeilberger approach, and the (fast) Zeilberger algorithm 
with which all above identities can be handled as well as 
a large list of identities that were studied by Gessel and Stanton (1982).

Our extensions therefore unify the verification of hypergeometric identities.

\section{The Gosper Algorithm}
\label{sec:The Gosper Algorithm}

In this section we recall the celebrated Gosper algorithm 
(Gosper, 1978), see also Graham, Knuth and Patashnik (1994).

The Gosper algorithm deals with the question to find
an antidifference $s_k$ for given $a_k$, i.\ e.\ a sequence $s_k$ for
which
\be
a_k=s_k-s_{k-1}
\label{eq:G1}
\ee
in the particular case that $s_k$ is a hypergeometric term,
therefore
\be
\frac{s_k}{s_{k-1}}
\quad\mbox{is a rational function with respect to } k\;,
\label{eq:G2}
\ee
i.\ e.\ $s_{k}/s_{k-1}\in\Q(k)$. We call this {\sl indefinite summation}.

Note that if a hypergeometric term antidifference $s_k$ exists,
we call the input function $a_k$ {\sl Gosper-summable} which then
itself is a hypergeometric term since by (\ref{eq:G1}) and (\ref{eq:G2})
\[
\frac{a_k}{a_{k-1}}=\frac{s_k-s_{k-1}}{s_{k-1}-s_{k-2}}=
\frac{\frac{s_k}{s_{k-1}}-1}{1-\frac{s_{k-2}}{s_{k-1}}}
=\frac{u_k}{v_k}
\in\Q(k)
\]
is rational, i.\ e.\ $u_k, v_k\in\Q[k]$ are polynomials.

Now, Gosper uses a representation lemma for rational functions to
express $a_k/a_{k-1}$ in terms of polynomials.

The idea behind this step comes from the following observation:
If we calculate $a_k$ from $s_k=(2k)!/k!$, e.g., we get 
\[
a_k=s_k-s_{k-1}=\frac{(2k)!}{k!}-\frac{(2k-2)!}{(k-1)!}
=((2k)(2k-1)-k)\cdot\frac{(2k-2)!}{k!}
=k(4k-3)\cdot\frac{(2k-2)!}{k!}
,
\]
i.\ e.\ a product of a polynomial $p_k=k(4k-3)$ and a factorial term
$b_k=\frac{(2k-2)!}{k!}$ for which $b_k/b_{k-1}=q_k/r_k$ is rational,
and therefore $q_k$ and $r_k$ can be assumed to be polynomials.

Gosper shows then that such a representation with the property
\be
\gcd\:(q_k,r_{k+j})=1
\quad\quad\mbox{for all}\;j\in\N_0
\label{eq:gcd1}
\ee
generally can be found and gives an algorithm to generate it.

Therefore we have for $a_k$ the relation
\be
\frac{a_k}{a_{k-1}}=\frac{p_k}{p_{k-1}}\,\frac{q_k}{r_k}
\label{eq:G5}
\;,
\ee
$p_k$ corresponding to the polynomial and
$(q_k, r_k)$ to the factorial part of $a_k$.
 
Gosper finally defines the function $f_k$ by the equation
\be
s_k=\frac{q_{k+1}}{p_k}\,f_k\,a_k
\label{eq:G7}
\ee
for which one sees immediately that
\[
f_k=\frac{p_k}{q_{k+1}}\frac{s_k}{a_k}=
\frac{p_k}{q_{k+1}}\frac{s_k}{s_k-s_{k-1}}=
\frac{p_k}{q_{k+1}}\frac{\frac{s_k}{s_{k-1}}}{\frac{s_k}{s_{k-1}}-1}
\]
is rational. Using (\ref{eq:gcd1}),
Gosper proves the essential fact that $f_k$ is a polynomial.

It follows from its defining equation that the polynomial $f_k$ 
satisfies 
\[
a_k=s_k-s_{k-1}=\frac{q_{k+1}}{p_k}\,f_k\,a_k-
\frac{q_{k}}{p_{k-1}}\,f_{k-1}\,a_{k-1}
\;,
\]
or multiplying by $p_k/a_k$, and using (\ref{eq:G5}), one gets
the recurrence equation
\be
p_k=q_{k+1}\,f_k-q_{k}\,\frac{p_k}{p_{k-1}}\frac{a_{k-1}}{a_k}\,f_{k-1}
=
q_{k+1}\,f_k-r_k\,f_{k-1}
\;.
\label{eq:fequation}
\ee
Using (\ref{eq:fequation}),
Gosper gives an upper bound for the degree of $f$ in terms
of the degrees of $p_k$, $q_k$, and $r_k$ which yields a fast method
to calculate $f_k$, so that we finally find $s_k$, given by (\ref{eq:G7}).
This shows in particular that whenever
$a_k$ possesses a closed form antidifference $s_k$ then necessarily
$s_k$ is a rational multiple of $a_k$:
\[
s_k=R_k\,a_k
\quad\quad\mbox{with}\quad\quad
R_k=\frac{q_{k+1}\,f_k}{p_k}
\;.
\]
If any of the steps to find the polynomial $f_k$ fails, the algorithm proves 
that no hypergeometric term antidifference $s_k$ of $a_k$ exists.

Therefore the Gosper algorithm is a {\sl decision procedure} which either
returns ``No closed form antidifference exists'' or returns a closed 
form antidifference $s_k$ of $a_k$, {\sl provided one can decide
the rationality of $a_k/a_{k-1}$}, i.\ e.\ one finds polynomials $u_k, v_k$
such that $a_k/a_{k-1}=u_k/v_k$. 
In so far, the Gosper algorithm is an algorithm with input $u_k$ and $v_k$
rather than $a_k$. 

Since without preprocessing,
the user's input is $a_k$ rather than the polynomials $u_k$ and $v_k$,
the success of an implementation depends heavily on an algorithm 
quickly and safely calculating $(u_k, v_k)$ given $a_k$.
In Algorithm~\ref{algo:3}, we present such a method.
It turns out that none of the existing implementations
of Gosper's algorithm uses such a method, examples of which we will consider
later.

In case, the Gosper algorithm provides us with an antidifference 
$s_k$ of $a_k$, any sum
\[
\sum_{k=m}^{n} a_k=s_{n}-s_{m-1}
\]
can be easily calculated by an evaluation of $s_k$ at the boundary points
like in the integration case. Note, however, that the sum
\begin{equation}
\sum_{k=0}^n \ueber{n}{k}
\label{eq:nchoosek}
\end{equation}
e.\ g.\
is not of this type as the summand $\ueber{n}{k}$ depends on the upper
boundary point $n$ explicitly. This is an example of a definite sum
that we consider in \S~\ref{sec:The Wilf-Zeilberger Method}.

Gosper implemented his algorithm in the \Macsyma\ {\tt nusum} command,
an implementation of the algorithm is distributed
with the {\tt sum} command of the \Maple\ system (to check its use set
{\tt infolevel[sum]:=5}), and one was delivered with
\Mathematica\ Version 1.2 (in the package
{\tt Algebra/GosperSum.m}). Another 
\Mathematica\ implementation was given by Paule and Schorn (1994).

On the lines of Koornwinder (1993), together with Gregor St\"olting
I implemented the Gosper algorithm
in \Reduce\ (Koepf, 1994) and \Maple,
using the simple decision procedure for rationality of
hypergeometric terms described in Algorithm \ref{algo:3} below rather
than internal simplification procedures (like \Maple's {\tt expand}). 
In \S~\ref{sec:The Wilf-Zeilberger Method}, we will show that
this makes our implementations that can be obtained from the author 
much stronger than the previous ones.

It is almost trivial but decisive that the following is a decision procedure
for the rationality of $a_k/a_{k-1}$ for input $a_k$
(at least) of a special type:
\begin{algorithm}
\label{algo:3}
{\rm
({\tt simplify\verb+_+combinatorial})\\
The following algorithm decides the rationality of $a_k/a_{k-1}$:
\begin{enumerate}
\item
Input: $a_k$ as ratio of products of
rational functions, exponentials, factorials, $\Gamma$ function terms, binomial
coefficients, and Pochhammer symbols that are rational-linear in their
arguments.
\item
({\tt togamma})\\
Build $a_k/a_{k-1}$, and convert all occurrences of factorials,
binomial coefficients, and Poch\-hammer symbols to $\Gamma$ function terms.
\item
({\tt simplify\verb+_+gamma})\\
Recursively rewrite this expression according to the rule
\[
\Gamma\:(a+k)=(a)_k\cdot \Gamma\:(a)
\]
($(a)_k:=a(a+1)\cdots(a+k-1)$ denoting the Pochhammer symbol)
whenever the arguments $a$ and $a+k$ of two representing $\Gamma$ function
terms have positive integer difference $k$. Reduce the final fraction
cancelling common $\Gamma$ terms.
\item
({\tt simplify\verb+_+power})\\
Recursively rewrite the last expression according to the rule
\[
b^{a+k}=b^k\,b^a
\]
whenever the arguments $a$ and $a+k$ of two representing exponential terms
have positive integer difference $k$. Reduce the final fraction
cancelling common exponential terms.
\item
The expression $a_k/a_{k-1}$ is rational if and only if the resulting
expression in step 4 is rational $u_k/v_k$, $u_k,v_k\in\Q[k]$.
\item
Output: $(u_k, v_k)$.
\end{enumerate}
}
\end{algorithm}
Note that this result follows immediately from the given form of 
$a_k$ and therefore of the
expression $a_k/a_{k-1}$ considered. 

As an example, the rationality of $a_k/a_{k-1}$ of
\[
a_k=\frac{\Gamma\:(2 k)}{4^k\,\Gamma\:(k)\,\Gamma\:(k+1/2)}
\]
is recognized by the given procedure, and from the resulting information
($a_k/a_{k-1}=1$), by induction $a_k=1/(2\sqrt\pi)$ 
(Abramowitz and Stegun, 1964, (6.1.18)).

Algorithm \ref{algo:3} does also apply to
\[
a_k=\Gamma\:(2 k)-\alpha\,4^k\,\Gamma\:(k)\,\Gamma\:(k+1/2)
\;,
\]
and leads to
\[
\frac{a_k}{a_{k-1}}=
2\,\left (2\,k-1\right )\left (k-1\right )
\]
which is true whenever $\alpha\neq \frac{1}{2\sqrt\pi}$.
If $\alpha= \frac{1}{2\sqrt\pi}$, however,
$a_k\equiv 0$, and therefore $a_k/a_{k-1}$ is not
properly defined.

For this reason, Gosper's algorithm with input 
\[
a_k=\Gamma\:(2 k)-\frac{1}{2\sqrt\pi}\,4^k\,\Gamma\:(k)\,\Gamma\:(k+1/2)
\]
fails to find the true antidifference $s_k=1$.

We note, however,
that in most cases also sums of ratios of the described form can be treated
by the same method without using multiplication formulas of the $\Gamma$
function explicitly.
An important family of examples of this type will be considered in 
the next section.

We note, finally,
that an implementation of Algorithm \ref{algo:3}
in a computer algebra system allows the user to enter his input in the
form in which it is found in the literature, and no preprocessing is
necessary. The user can be sure that the rationality is decided
correctly. Unfortunately, this is not so with any of
the current implementations.

\Maple's {\tt expand} command which is used for that purpose 
in Zeilberger's (1990) and Koornwinder's (1993) implementations of
Zeilberger's algorithm (see \S~\ref{sec:The Zeilberger Algorithm}), e.\ g., does
not the job required. The same is valid for \Mathematica's {\tt
FactorialSimplify} procedure
that comes in the package {\tt DiscreteMath/RSolve.m}.
Also, the Gosper implementation which
comes with \Maple's {\tt sum} command, has the same failure. Paule and
Schorn's \Mathematica\ implementation aborts in those cases with
the message {\tt input not interpretable}, whereas Gosper's {\tt nusum}
command gives the (wrong)
error message {\tt errexp1 NON-RATIONAL TERM RATIO TO NUSUM}.

Therefore in all these implementations, the fact that Gosper's algorithm is a 
decision procedure, unfortunately is completely lost. 
An example for that fact is given by the expression
\[
a_k=\frac{1}{2^n}\ueber{n}{k}-\frac{1}{2^{n-1}}\ueber{n-1}{k}
\]
with
\[
\frac{a_k}{a_{k-1}}=
{\frac {\left (n-k+1\right )\left (n-2\,k\right )}{k\left (n-2\,k+2
\right )}}
\;,
\]
for which our Gosper implementations succeed very quickly,
whereas \Maple's {\tt sum} command as well as \Mathematica's {\tt
GosperSum}, Paule/Schorn's {\tt Gosper}, and Gosper's {\tt nusum}
fail if the input is not preprocessed by the user.

This is a very simple example of an important type of examples considered 
next.

\section{The Wilf-Zeilberger Method}
\label{sec:The Wilf-Zeilberger Method}

Examples for an application of the Gosper algorithm in connection
with Algorithm \ref{algo:3} are given by the
Wilf-Zeilberger method on {\sl definite summation} (Wilf and Zeilberger, 1990),
see also Wilf (1993).

The Wilf-Zeilberger method is a direct application of Gosper's algorithm
to prove identities of the form
\be
s_n:=\sum_{k\in\ZZ} F(n,k)=1
\label {eq:WZ1}
\ee
for which $F(n,k)$ is a hypergeometric term with respect to both $n$ and $k$,
i.\ e.\
\[
\frac{F(n,k)}{F(n-1,k)}
\quad
\mbox{and}
\quad
\frac{F(n,k)}{F(n,k-1)}
\quad
\mbox{are rational functions with respect to both $n$ and $k$,}
\]
where $n$ is assumed to be an integer, and the sum is to be taken over
all integers $k\in\Z$.

To prove a statement of the form (\ref{eq:WZ1})
by the WZ method%
\footnote{Note that Wilf and Zeilberger use forward differences
rather that downward differences, whereas we decided to follow
Gosper's original treatment. There is no theoretical difference
between these two approaches, though.}, 
one applies Gosper's algorithm to the expression 
\[
a_k:=F(n,k)-F(n-1,k)
\]

{\small
\begin{table}[hbtp]
\caption{The WZ method}
\label{table:WZ}
\vspace*{3mm}
\begin{tabular}{lrp{12cm}}
Theorem&$n$&$R(n,k)$\\[2mm]
\hline \\[-2mm]
Vandermonde & $-a$ &
$\displaystyle
-{\frac {\left (b+k\right )\left (-n+k\right )}{n\left (c+n-1\right )}}
$
\\[4mm]
\hline \\[-2mm]
Saalsch\"utz & $n$ &
$\displaystyle
-{\frac {\left (b+k\right )\left (-n+k\right )\left (a+k\right )}{n
\left (c+n-1\right )\left (1+a+b-c-n+k\right )}}
$
\\[4mm]
\hline \\[-2mm]
Kummer & $-b$ &
$\displaystyle
{\frac {\left (a+k\right )\left (-n+k\right )}{n\left (a+2\,n\right )}}
$
\\[4mm]
\hline \\[-2mm]
Dixon & $-c$ &
$\displaystyle
-{\frac {\left (a+k\right )\left (-n+k\right )\left (b+k\right )}{n
\left (a-b+n\right )\left (a+2\,n\right )}}
$
\\[4mm]
\hline \\[-2mm]
\multicolumn{1}{l}
{\hspace*{-2mm}\begin{tabular}{l}Watson\\ Whipple\end{tabular}}
& $-c$ &
$\displaystyle
{2\frac {\left (a+k\right )\left (-n+k\right )\left (b+k\right )}
{\left (-1+a+b+2\,n\right )\left (-2\,n+1+k\right )\left (-2\,n+k
\right )}}
$
\\[4mm]
\hline \\[-2mm]
Whipple & $-c$ &
$\displaystyle
-{\frac {\left (a+k\right )\left (a-1-k\right )\left (-n+k\right )}{n
\left (2-2\,n-e+k\right )\left (1-2\,n-e+k\right )}}
$
\\[4mm]
\hline \\[-2mm]
Dougall & $n$ &
$\displaystyle
{\frac {\left (2\,a-b-c-d+2\,n\right )\left (a+k\right )\left (-n+k
\right )\left (b+k\right )\left (c+k\right )\left (d+k\right )}{n
\left (a+2\,k\right )\left (a-b-c+n-d-k\right )\left (a-d+n\right )
\left (a-c+n\right )\left (a-b+n\right )}}
$
\\[4mm]
\hline \\[-2mm]
Dougall & $-e$ &
$\displaystyle
-{\frac {\left (a+k\right )\left (-n+k\right )\left (c+k\right )\left
(d+k\right )}{n\left (a+2\,k\right )\left (a-c+n\right )\left (a-d+n
\right )}}
$
\\[4mm]
\hline \\[-2mm]
Whipple & $-e$ &
$\displaystyle
{\frac {\left (d+k\right )\left (-n+k\right )\left (a+k\right )}{n
\left (a+2\,k\right )\left (a-d+n\right )}}
$
\\[4mm]
\hline \\[-2mm]
Bailey & $n$ &
$\displaystyle
-{\frac {\left ({a}^{2}\!+\!2\,a\!-\!wa\!+\!na\!+\!2\!-\!2\,w\!-
\!2\,kw\!+\!2\,ka\!+\!2\,k\!+\!2\,kn
\right )\left (a+k\right )\left (-n+k\right )}{\left (-w+a+n\right )n
\left (a+2\,k\right )\left (w+n-1\right )}}
$
\\[4mm]
\hline \\[-2mm]
Bailey & $n$ &
$\displaystyle
-{\frac {\left (-2\,b-2\,{b}^{2}+2\,nb+ab-1+n-k\right )\left (a+k
\right )\left (-n+k\right )\left (b+k\right )}{nb\left (1+2\,b-n+k
\right )\left (a-2\,b+2\,n-2\right )\left (a-b+n\right )}}
$
\\[4mm]
\hline \\[-2mm]
Bailey & $n$ &
$\displaystyle
-{\frac {\left (2\,b+ab+1-n+2\,kb+k\right )\left (b+k\right )\left (-n
+k\right )\left (a+k\right )}{nb\left (a+2\,k\right )\left (1+2\,b-n+k
\right )\left (a-b+n\right )}}
$
\\[4mm]
\hline \\[-2mm]
Bailey & $n$ &
$\displaystyle
-{\frac {
\left (a+k\right )\left (-n+k\right )\left (b+k\right )}{nb
\left (a+2\,k\right )\left (2+2\,b-n+k\right )\left (a-2\,b-3+2\,n
\right )\left (a-b+n\right )}}
\cdot
$
\\[3mm]
\multicolumn{3}{r}{
$\displaystyle
\cdot
\left (-8b\!-\!4{b}^{2}\!+\!6nb\!-\!ab\!-\!2{n}^{2}\!+\!2nba\!-\!4\!+\!
6n\!-\!2
{b}^{2}a\!+\!{a}^{2}b\!-\!6k\!-\!8kb\!-\!4{b}^{2}k\!+\!4kn\!+\!4kbn\!+\!2kba\!-
\!2{k}^{2
}\right )
$
}
\end{tabular}
\end{table}
}\noindent
with respect to the variable $k$. If successful, this generates $G(n,k)$ with
\be
a_k=F(n,k)-F(n-1,k)=G(n,k)-G(n,k-1)
\;,
\label{eq:WZ2}
\ee
and summing over all $k$ leads to
\[
s_n-s_{n-1}=\sum_{k\in\ZZ} \Big( F(n,k)-F(n-1,k)\Big)
=\sum_{k\in\ZZ} \Big( G(n,k)-G(n,k-1)\Big)=0
\]
since the right hand side is telescoping. Therefore $s_n$ is constant,
$s_n=s_0$, and if we are able to prove $s_0=1$, we are done.
Note that $s_0=1$ generally can be proved if the series considered
is terminating in which case also no questions concering convergence arise.

Since the WZ method only works if $n$ is an integer, we can try to
prove the statements of Bailey's list in Table \ref{table:Bailey}
anyway only if one of the upper parameters of the hypergeometric series
involved is a negative integer. The extension to the general case
is over the capabilities of the methods of this article,
and must be handled by other means.

Note that the rationality of $a_k/a_{k-1}$
for the WZ method is decided by Algorithm \ref {algo:3} since
\[
\frac{a_{k}}{a_{k-1}}=\frac{F(n,k)-F(n-1,k)}{F(n,k-1)-F(n-1,k-1)}
=\frac{F(n,k)}{F(n,k-1)}\cdot
\frac{1-\frac{F(n-1,k)}{F(n,k)}}{1-\frac{F(n-1,k-1)}{F(n,k-1)}}
\;.
\]

Note moreover, that the application of Gosper's algorithm may be slow.
But as soon as we have found the function $G(n,k)$, we easily can
calculate the rational function
\[
R(n,k):=\frac{G(n,k)}{F(n,k)}
\]
by an application of Algorithm \ref{algo:3}. $R(n,k)$ is rational
since the proof of Gosper's algorithm shows that $G(n,k)$ is a
rational multiple of $a_{k}=F(n,k)-F(n-1,k)$,
$G(n,k)=r(n,k)\cdot(F(n,k)-F(n-1,k))$, say, so that
\[
R(n,k)=\frac{G(n,k)}{F(n,k)}=
r(n,k) \frac{F(n,k)-F(n-1,k)}{F(n,k)}=
r(n,k)\lk 1-\frac{F(n-1,k)}{F(n,k)}\rk
\]
is rational. $R(n,k)$ is called the
{\sl rational certificate} of $F(n,k)$. Once the rational certificate
of a hypergeometric expression $F(n,k)$ is known, it is a matter of
pure rational arithmetic (which is fast) to decide the validity of
(\ref{eq:WZ1}) since the only thing that one has to show is
(\ref{eq:WZ2}) which after division by $F(n,k)$
is equivalent 
(modulo an application of Algorithm \ref{algo:3})
to the purely rational identity
\be
1-R(n,k)+R(n,k-1)\frac{F(n,k-1)}{F(n,k)}-\frac{F(n-1,k)}{F(n,k)}=0
\label{eq:rationalcertificate}
\;.
\ee
As an example, to prove the Binomial Theorem (compare (\ref{eq:nchoosek}))
in the form
\be
s_n:=\sum_{k=0}^n F(n,k)=\sum_{k=0}^n \frac{1}{2^n}\,\ueber{n}{k}=1
\label{eq:BinomialWZ}
\ee
by the WZ method, Algorithm~\ref{algo:3} yields
\[
\frac{a_k}{a_{k-1}}=\frac{F(n,k)-F(n-1,k)}{F(n,k-1)-F(n-1,k-1)}=
{\frac {\left (n-k+1\right )\left (n-2\,k\right )}{k\left (n-2\,k+2
\right )}}
\]
so that Gosper's algorithm can be applied, and results in
\[
G(n,k)=\frac{k}{2^n(2k - n)}\lk 2\ueber{n - 1}{k} - \ueber{n}{k}\rk
\;.
\]
This proves (\ref{eq:BinomialWZ}) since $s_0=\sum\limits_{k=0}^0 1=1$.

The rational certificate function is
\[
R(n,k)=
\frac{G(n,k)}{F(n,k)}=
\frac{k-n}{n}
\;,
\]
and the verification of identity (\ref{eq:BinomialWZ}) is therefore reduced to
simplify the rational expression
\[
1\!-\!R(n,k)\!+\!R(n,k-1)\frac{F(n,k-1)}{F(n,k)}\!-\!\frac{F(n-1,k)}{F(n,k)}=
1\!-\!\frac{k-n}{n}\!+\!\frac{k\!-\!1\!-\!n}{n}
\frac{k}{n\! +\! 1\!-\!k}\!-\!\frac{2(n\!-\!k)}{n}
\]
to zero.

Table \ref {table:WZ} is a complete list of those identities of Bailey's list
(Table \ref{table:Bailey}) that can be treated by the given method together
with their rational certificates with which the reader may verify them
easily.

Note that neither the statements of Gau{\ss} and Bailey of argument
$x=1/2$ (p.\ 11) are accessible with respect to any of the parameters involved,
nor can Watson's Theorem \mbox{(p.\ 16)} 
be proved by the WZ method with respect to
Watson's original integer parameter $a$, nor can the method be applied
to Whipple's Theorem (p.\ 16) concerning parameters $a$ or $b$ since in all
these cases the term ratio $a_k/a_{k-1}$ is not rational.

Note further that both our \Reduce\ and our 
\Maple\ implementations generate the results of 
Table~\ref {table:WZ}, and only the calculation of the rational certificate
of Dougall's Theorem needs more than a few seconds. In the appendix,
we will present some example results. On the other hand,
\Maple's {\tt sum} command was not successful with a single example
without preprocessing the input (entered in factorial or $\Gamma$ notation), 
in most cases quickly responding with the incorrect statement
{\tt Gosper's algorithm fails}, and unsuccessfully
trying other methods afterwards.

In \S~\ref{sec:An Extended Version of the Gosper Algorithm},
we consider a generalization of the WZ method. 
To be able to consider the most general case, 
we present an extended version of Gosper's algorithm next.

\section{An Extended Version of the Gosper Algorithm}
\label{sec:An Extended Version of the Gosper Algorithm}

Here we deal with the question, given a nonnegative integer $m$,
to find a sequence $s_k$ for given $a_k$ satisfying
\be
a_k=s_k-s_{k-m}
\label{eq:1}
\ee
in the particular case that $s_k$ is an {\sl $m$-fold hypergeometric term},
i.\ e.
\be
\frac{s_k}{s_{k-m}}
\quad\mbox{is a rational function with respect to } k\;.
\label{eq:2}
\ee
Note that in the given case the input function $a_k$ itself
is an $m$-fold hypergeometric term since by (\ref{eq:1}) and (\ref{eq:2})
\[
\frac{a_k}{a_{k-m}}=\frac{s_k-s_{k-m}}{s_{k-m}-s_{k-2m}}=
\frac{\frac{s_k}{s_{k-m}}-1}{1-\frac{s_{k-2m}}{s_{k-m}}}
=\frac{u_k}{v_k}
\label{eq:initial choice}
\]
is rational, i.\ e., $u_k$ and $v_k$ can be chosen to be polynomials.

Assume first, given $a_k$, we have found $s_k$ with $s_k-s_{k-m}=a_k$. Then
we can easily construct an antidifference $\tilde s_k$ of $a_k$ by
\be
\tilde s_k:=s_k+s_{k-1}+\cdots+s_{k-(m-1)}
\label{eq:antifrom-m}
\ee
since then
\[
\tilde s_k-\tilde s_{k-1}=
(s_k+\cdots+s_{k-(m-1)})-
(s_{k-1}+\cdots+s_{k-m})=
s_k-s_{k-m}=a_k
\;.
\]
Note, however, that in general, this antidifference is not a 
hypergeometric term, but is a finite sum of hypergeometric terms.

Assume next that given $a_k$ there exists a hypergeometric term $s_k$ with
$s_k-s_{k-m}=a_k$. Then obviously also
\[
\frac{s_k}{s_{k-m}}=
\frac{s_k}{s_{k-1}}\cdot\frac{s_{k-1}}{s_{k-2}}\cdots\frac{s_{k-(m-1)}}{s_{k-m}}
\]
is rational, and therefore our algorithm below will find $s_k$.

An $m$-fold antidifference always can be constructed by
an application of Gosper's original algorithm in the
following way:
\begin{algorithm}
\label{algo:1}
{\rm
({\tt extended\verb+_+gosper})\\[3mm]
The following steps generate an $m$-fold antidifference:
\begin{enumerate}
\item
Input: $a_k$, and $m\in\N$.
\item
Define $b_k:=a_{km}$.
\item
Apply Gosper's algorithm to $b_k$ with respect to $k$. Get the antidifference
$t_k$ of $b_k$, or the statement: ``No hypergeometric term antidifference of 
$b_k$, and therefore no $m$-fold hypergeometric term antidifference of $a_k$ 
exists.''
\item
The output $s_k:=t_{k/m}$ is a solution of (\ref{eq:1}) with the property
(\ref{eq:2}).
\end{enumerate}
}
\end{algorithm}
\begin{proof}
This is valid since $a_k=s_k-s_{k-m}$ is equivalent to
\[
b_k=a_{km}=s_{km}-s_{km-m}=
t_k-t_{k-1}
\]
and since
\[
\frac{t_{k}}{t_{k-1}}=\frac{s_{km}}{s_{km-m}}
\]
describes the transformation between $t_k$ and $s_k$.
\end{proof}
\\[-3mm]
As an example, we consider $a_k:=k\lk\frac{k}{2}\rk!$, and $m=2$. 
Then $b_k=a_{2k}= 2k\,k!$, and Gosper's algorithm yields $t_k=2(k+1)k!$. 
Therefore $s_k=t_{k/2}=(k+2)\lk\frac{k}{2}\rk!$ has the property that
\[
s_k-s_{k-2}=a_k
\;.
\]
By (\ref{eq:antifrom-m}), we moreover find the antidifference
\[
\tilde s_k=s_k+s_{k-1}=
(k+2)\lk\frac{k}{2}\rk!+(k+1)\lk\frac{k-1}{2}\rk!
\]
of $a_k$.

We consider two other examples: If $a_k=\ueber{k/3}{n}$ then
our algorithm generates the antidifference
\[
s_k=\frac{1}{3(n+1)}\lk (k+3) \ueber{\frac{k}{3}}{n} + 
(k+2)\ueber{\frac{k-1}{3}}{n}+ (k+1)\ueber{\frac{k-2}{3}}{n} \rk
\;,
\]
and if $a_k=\ueber{n}{k/2}-\ueber{n}{k/2-1}$ then
\begin{eqnarray*}
s_k&=&
\frac{(2n + 3 - k)(n + 1 - k) }{2(n + 2 - k)(n + 1 - k)}
\lk \ueber{n}{\frac{k-1}{2}}-\ueber{n}{\frac{k-3}{2}}\rk
\\
&&+ \frac{(n + 2 - k)(2n + 2 - k)}{2(n + 2 - k)(n + 1 - k)}\lk
\ueber{n}{k/2}-\ueber{n}{k/2-1}\rk
\;.
\end{eqnarray*}
Note, however, that we will use $m$-fold hypergeometric antidifferences
rather than non-hyper\-geometric antidifferences in the later chapters. 

Now, we give an algorithm that finds an appropriate nonnegative integer $m$
for an arbitrary input function $a_k$ given as ratio of products of
rational functions, exponentials, factorials, $\Gamma$ function terms, binomial
coefficients, and Pochhammer symbols that are rational-linear in their
arguments:
\begin{algorithm}
\label{algo:2}
{\rm
({\tt find\verb+_+mfold})\\[3mm]
The following is an algorithm generating a successful choice for $m$
for an application of Algorithm~\ref{algo:1}.
\begin{enumerate}
\item
Input: $a_k$ as ratio of products of
rational functions, exponentials, factorials, $\Gamma$ function terms, binomial
coefficients, and Pochhammer symbols that are rational-linear in their
arguments.
\item
Build the list of all arguments. They are of the form
$p_j/q_j\,k+\al_j$ with integer $p_j$ and $q_j$, $p_j/q_j$ in
lowest terms, $q_j$ positive.
\item
Calculate $m:=\lcm\{q_j\}$.
\end{enumerate}
}
\end{algorithm}

\begin{proof}
It is clear that the procedure generates a representation for $b_k=a_{km}$
with the given choice of $m$ which is integer-linear in the arguments involved.
Since in this case $b_k/b_{k-1}$ is rational,
Algorithm~\ref{algo:1} is applicable.
\end{proof}
\\[-3mm]
We mention that in our example cases above, the given procedure yields the
desired values $m=2$ for $a_k:=k\lk\frac{k}{2}\rk!$, $m=3$ for
$a_k=\ueber{k/3}{n}$, and $m=2$ for 
$a_k=\ueber{n}{k/2}-\ueber{n}{k/2-1}$.

\section{Extension of the WZ method}
\label{sec:Extension of the WZ method}

In this section we will give an extended version of the WZ method which
resolves the questions that remained open in 
\S~\ref{sec:The Wilf-Zeilberger Method}
so that finally Bailey's complete list (Table
\ref{table:Bailey}) can be settled using a unifying approach.

Assume that for a hypergeometric identity the WZ method fails. This may happen
either because $a_k/a_{k-1}$ is not rational, or because there is no single 
formula for the result like in Andrews' statement
\be
_3 F_2\left.
\!\!
\left(
\!\!\!\!
\begin{array}{c}
\multicolumn{1}{c}{\begin{array}{c}
-n\;, n+3a\;, a
\end{array}}\\[1mm]
\multicolumn{1}{c}{\begin{array}{c}
3a/2\;,(3a+1)/2
            \end{array}}\end{array}
\!\!\!\!
\right| \frac{3}{4}\right)
=
\funkdef{0}{n\neq 0 {\mbox{ (mod }} 3)}{\displaystyle
\frac{n!\,(a+1)_{n/3}}
{(n/3)!\,(3a+1)_n}
}
\label{eq:Gessel-Stanton (1.1)}
\ee
which---together with many similar statements---can be found in a paper
of Gessel and Stanton (1982), Equation (1.1).

In such cases, we proceed as follows: To prove an identity of the form
\be
s_n:=\sum_{k\in\ZZ} F(n,k)={\rm constant}
\quad\quad{(n \mbox{ mod }} m \mbox{ constant})
\label {eq:WZ1extended}
\;,
\ee
$m$ denoting a certain positive integer,
$F(n,k)$ being an {\sl $(m,l)$-fold hypergeometric term} with respect to
$(n,k)$, i.\ e.\
\[
\frac{F(n,k)}{F(n-m,k)}
\quad
\mbox{and}
\quad
\frac{F(n,k)}{F(n,k-l)}
\quad
\mbox{are rational functions with respect to both $n$ and $k$,}
\]
and $n$ assuming to be an integer,
we apply our extended version of Gosper's algorithm to find an
$l$-fold antidifference of the expression
\[
a_k:=F(n,k)-F(n-m,k)
\]

{\small
\begin{table}[hbtp]
\caption{Gessel and Stanton's hypergeometric identities}
\label{table:Gessel-Stanton}
\vspace*{3mm}
\begin{tabular}{rlp{14.5cm}}
Eq. & Identity\\[2mm] \hline \\[-2mm]
(1.1) 
&
$\displaystyle
_3 F_2\left.
\!\!
\left(
\!\!\!\!
\begin{array}{c}
\multicolumn{1}{c}{\begin{array}{c}
-n\;, n+3a\;, a
\end{array}}\\[1mm]
\multicolumn{1}{c}{\begin{array}{c}
3a/2\;,(3a+1)/2
            \end{array}}\end{array}
\!\!\!\!
\right| \frac{3}{4}\right)
=
\funkdef{0}{n\neq 0 {\mbox{ (mod }} 3)}{\displaystyle
\frac{n!\,(a+1)_{n/3}}
{(n/3)!\,(3a+1)_n}
}
$
\\[5mm]
(1.2)
&
$\displaystyle
_5 F_4\left.
\!\!
\left(
\!\!\!\!
\begin{array}{c}
\multicolumn{1}{c}{\begin{array}{c}
2a\;, 2b\;, 1-2b\;,1+2a/3\;, -n
\end{array}}\\[1mm]
\multicolumn{1}{c}{\begin{array}{c}
a-b+1\;, a+b+1/2\;,2a/3\;, 1+2a+2n
            \end{array}}\end{array}
\!\!\!\!
\right| \frac{1}{4}\right)
=
\frac{(a+1/2)_{n}\,(a+1)_n}
{(a+b+1/2)_n\,(a-b+1)_n}
$
\\[5mm]
(1.3)
&
$\displaystyle
_5 F_4\left.
\!\!
\left(
\!\!\!\!
\begin{array}{c}
\multicolumn{1}{c}{\begin{array}{c}
a, b, a\!+\!1/2\!-\!b,1\!+\!2a/3, -n
\end{array}}\\[1mm]
\multicolumn{1}{c}{\begin{array}{c}
2a\!+\!1\!-\!2b,2b,2a/3,1\!+\!a\!+\!n/2
            \end{array}}\end{array}
\!\!\!\!
\right| 4\right)
=
\funkdef{0}{n\mbox{ odd}}{\displaystyle
\frac{n!\,(a+1)_{n/2}\,2^{-n}}
{(\frac{n}{2})!(a\!-\!b\!+\!1)_{n/2}(b\!+\!\frac{1}{2})_{n/2}}
}
$
\\[5mm]
(1.4)
&
$\displaystyle
_3 F_2\left.
\!\!
\left(
\!\!\!\!
\begin{array}{c}
\multicolumn{1}{c}{\begin{array}{c}
1/2+3a\;, 1/2-3a\;, -n
\end{array}}\\[1mm]
\multicolumn{1}{c}{\begin{array}{c}
1/2\;, -3n
            \end{array}}\end{array}
\!\!\!\!
\right| \frac{3}{4}\right)
=
\frac{(1/2-a)_{n}\,(1/2+a)_n}
{(1/3)_n\,(2/3)_n}
$
\\[5mm]
(1.5)
&
$\displaystyle
_3 F_2\left.
\!\!
\left(
\!\!\!\!
\begin{array}{c}
\multicolumn{1}{c}{\begin{array}{c}
1+3a\;, 1-3a\;, -n
\end{array}}\\[1mm]
\multicolumn{1}{c}{\begin{array}{c}
3/2\;, -1-3n
            \end{array}}\end{array}
\!\!\!\!
\right| \frac{3}{4}\right)
=
\frac{(1+a)_{n}\,(1-a)_n}
{(2/3)_n\,(4/3)_n}
$
\\[5mm]
(1.6)
&
$\displaystyle
_3 F_2\left.
\!\!
\left(
\!\!\!\!
\begin{array}{c}
\multicolumn{1}{c}{\begin{array}{c}
2a\;, 1-a\;,-n
\end{array}}\\[1mm]
\multicolumn{1}{c}{\begin{array}{c}
2a+2\;, -a-1/2-3n/2
            \end{array}}\end{array}
\!\!\!\!
\right| 1\right)
=
\frac{((n+3)/2)_{n}\,(n+1)(2a+1)}
{(1+(n+2a+1)/2)_n\,(2a+n+1)}
$
\\[5mm]
(1.7)
&
$\displaystyle
_7 F_6\left.
\!\!
\left(
\!\!\!\!
\begin{array}{c}
\multicolumn{1}{c}{\begin{array}{c}
2a\;, 2b\;, 1-2b\;, 1+2a/3\;, a+d+n+1/2\;, a-d\;, -n
\end{array}}\\[1mm]
\multicolumn{1}{c}{\begin{array}{c}
a-b+1\;, a+b+1/2\;, 2a/3\;, -2d-2n\;, 2d+1\;, 1+2a+2n
            \end{array}}\end{array}
\!\!\!\!
\right| 1\right)
=$\\[4mm]
&$
\hspace*{0.5cm}
\displaystyle
\frac{(2a+1)_{2n}\,(b+d+1/2)_n\,(d-b+1)_n}
{(2d+1)_{2n}\,(a+b+1/2)_n\,(a-b+1)_n}
=
\frac{(a+1/2)_n\,(a+1)_n\,(b+d+1/2)_n\,(d-b+1)_n}
{(a+b+1/2)_n\,(a-b+1)_n\,(d+1/2)_n\,(d+1)_n}
$
\\[5mm]
(1.8)
&
$\displaystyle
_7 F_6\left.
\!\!
\left(
\!\!\!\!
\begin{array}{c}
\multicolumn{1}{c}{\begin{array}{c}
a\;, b\;, a+1/2-b\;, 1+2a/3\;, 1-2d\;,2a+2d+n\;, -n
\end{array}}\\[1mm]
\multicolumn{1}{c}{\begin{array}{c}
2a-2b+1\;, 2b\;, 2a/3\;, a+d+1/2\;, 1-d-n/2\;, 1+a+n/2
            \end{array}}\end{array}
\!\!\!\!
\right| 1\right)
=$\\[4mm]
&$
\hspace*{1.5cm}
\displaystyle
\funkdef{0}{n\;\mbox{odd}}{\displaystyle
\frac{(b+d)_{n/2}\,(d-b+a+1/2)_{n/2}\,n!\,(a+1)_{n/2}\,2^{-n}}
{(b+1/2)_{n/2}\,(a+d+1/2)_{n/2}\,(d)_{n/2}\,(n/2)!\,(a-b+1)_{n/2}}
}
$
\\[5mm]
(3.7)
&
$\displaystyle
_2 F_1\left.
\!\!
\left(
\!\!\!\!
\begin{array}{c}
\multicolumn{1}{c}{\begin{array}{c}
-n\;, -2n-2/3
\end{array}}\\[1mm]
\multicolumn{1}{c}{\begin{array}{c}
4/3
            \end{array}}\end{array}
\!\!\!\!
\right| -8\right)
=
\frac{(5/6)_{n}}{(3/2)_n}\left(-27\right)^n
$
\\[5mm]
(5.21)
&
$\displaystyle
_3 F_2\left.
\!\!
\left(
\!\!\!\!
\begin{array}{c}
\multicolumn{1}{c}{\begin{array}{c}
3a+1/2\;, 3a+1\;, -n
\end{array}}\\[1mm]
\multicolumn{1}{c}{\begin{array}{c}
6a+1\;, -n/3+2a+1
            \end{array}}\end{array}
\!\!\!\!
\right| \frac{4}{3}\right)
=
\funkdef{0}{n\neq 0 {\mbox{ (mod }} 3)}{\displaystyle
\frac{(1/3)_{n/3}\,(2/3)_{n/3}}
{(1+2a)_{n/3}\,(-2a)_{n/3}}
}
$
\\[5mm]
(5.22)
&
$\displaystyle
_2 F_1\left.
\!\!
\left(
\!\!\!\!
\begin{array}{c}
\multicolumn{1}{c}{\begin{array}{c}
-n\;, 1/2
\end{array}}\\[1mm]
\multicolumn{1}{c}{\begin{array}{c}
2n+3/2
            \end{array}}\end{array}
\!\!\!\!
\right| \frac{1}{4}\right)
=
\frac{(1/2)_{n}}{(2n+3/2)_n}\left(\frac{27}{4}\right)^n
$
\\[5mm]
(5.23)
&
$\displaystyle
_2 F_1\left.
\!\!
\left(
\!\!\!\!
\begin{array}{c}
\multicolumn{1}{c}{\begin{array}{c}
-n\;, -1/3-2n
\end{array}}\\[1mm]
\multicolumn{1}{c}{\begin{array}{c}
2/3
            \end{array}}\end{array}
\!\!\!\!
\right| -8\right)
=
(-27)^n
$
\\[5mm]
(5.24)
&
$\displaystyle
_2 F_1\left.
\!\!
\left(
\!\!\!\!
\begin{array}{c}
\multicolumn{1}{c}{\begin{array}{c}
-n\;, n/2+1
\end{array}}\\[1mm]
\multicolumn{1}{c}{\begin{array}{c}
4/3
            \end{array}}\end{array}
\!\!\!\!
\right| \frac{8}{9}\right)
=
\funkdef{0}{n\;\mbox{odd}}{\displaystyle
\frac{(1/2)_{n/2}}{(7/6)_{n/2}}\,(-3)^{-(n/2)}
}
$
\\[5mm]
(5.25)
&
$\displaystyle
_2 F_1\left.
\!\!
\left(
\!\!\!\!
\begin{array}{c}
\multicolumn{1}{c}{\begin{array}{c}
-n\;, 1/2
\end{array}}\\[1mm]
\multicolumn{1}{c}{\begin{array}{c}
(n+3)/2
            \end{array}}\end{array}
\!\!\!\!
\right| 4\right)
=
\funkdef{0}{n\;\mbox{odd}}{\displaystyle
\frac{(1/2)_{n/2}\,(3/2)_{n/2}}{(5/6)_{n/2}\,(7/6)_{n/2}}
}
$
\\[5mm]
(5.27)
&
$\displaystyle
_4 F_3\left.
\!\!
\left(
\!\!\!\!
\begin{array}{c}
\multicolumn{1}{c}{\begin{array}{c}
1/3-n\;, -n/2\;, (1-n)/2\;, 22/21-3n/7
\end{array}}\\[1mm]
\multicolumn{1}{c}{\begin{array}{c}
5/6\;, 4/3\;, 1/21-3n/7
            \end{array}}\end{array}
\!\!\!\!
\right| -27\right)
=
\frac{(-8)^n}{1-9n}
$
\end{tabular}
\end{table}
}\noindent

{\small
\begin{table}[hbtp]
\caption{The extended WZ method}
\label{table:WZextended}
\vspace*{3mm}
\begin{tabular}{lrrp{10cm}}
Bailey p.&$n$&$m$&$R(n,k)$\\[2mm]
\hline \\[-2mm]
11, Gau{\ss} & $-a$ & $2$ &
$\displaystyle
-{\frac {\left (b+k\right )\left (n-k\right )}{\left (-b+n-1-2\,k
\right )n}}
$
\\[3mm]
\hline \\[-3mm]
11, Bailey & $-a$ & $2$ &
$\displaystyle
{\frac {\left (2\,n-1\right )\left (n-k\right )}{\left (c+n-1\right )
\left (n+k\right )}}
$
\\[3mm]
\hline \\[-3mm]
16, Watson & $-a$ & $2$ &
$\displaystyle
-{2\frac {\left (c+k\right )\left (b+k\right )\left (n-k\right )}
{\left (-1+n+2\,c\right )\left (-b+n-1-2\,k\right )n}}
$
\\[3mm]
\hline \\[-3mm]
16,  Whipple & $-a$ & $2$ &
$\displaystyle
{2\frac {\left (2\,n-1\right )\left (n-k\right )\left (c+k\right )}{
\left (2\,c-e+n\right )\left (-1+n+e\right )\left (n+k\right )}}
$
\\[3mm]
\hline \\[-3mm]
G.-S. Eq.&$m$&$R(n,k)$
\\[2mm]
\hline \\[-3mm]
(1.1) & $3$ &
\multicolumn{2}{l}{
$\displaystyle
{3\frac {\left (a+k\right )\left (n-k\right )\left (3\,a+2\,n-3
\right )}{\left (n+3\,a+k-2\right )\left (n+3\,a+k-1\right )n}}
$
}
\\[3mm]
\hline \\[-3mm]
(1.2) & $1$ &
\multicolumn{2}{l}{
$\displaystyle
-{\frac {\left (2\,a+k\right )\left (n-k\right )\left (2\,b-1-k\right 
)\left (2\,b+k\right )}{2\,n\left (3\,k+2\,a\right )\left (2\,a+2\,b-1
+2\,n\right )\left (a-b+n\right )}}
$
}
\\[3mm]
\hline \\[-3mm]
(1.3) & $2$ &
\multicolumn{2}{l}{
$\displaystyle
{4\frac {\left (b+k\right )\left (a+k\right )\left (2\,a-2\,b+1+2
\,k\right )\left (n-k\right )}{n\left (3\,k+2\,a\right )\left (2\,b-1+
n\right )\left (2\,a-2\,b+n\right )}}
$
}
\\[3mm]
\hline \\[-3mm]
(1.4) & $1$ &
\multicolumn{2}{l}{
$\displaystyle
{3\frac {\left (n-k\right )\left (6\,a-1-2\,k\right )\left (6\,a+
1+2\,k\right )}{\left (12\,n-4\,k\right )\left (3\,n-1-k\right )\left 
(3\,n-2-k\right )}}
$
}
\\[3mm]
\hline \\[-3mm]
(1.5) & $1$ &
\multicolumn{2}{l}{
$\displaystyle
{3\frac {\left (n-k\right )\left (3\,a-k-1\right )\left (3\,a+1+k
\right )}{\left (3\,n-1-k\right )\left (3\,n-k\right )\left (3\,n-k+1
\right )}}
$
}
\\[3mm]
\hline \\[-3mm]
(1.6) & $2$ &
\multicolumn{2}{l}{
$\displaystyle
{\frac {\left (-4\,a\!+\!4\,an\!+\!18\,{n}^{2}\!-\!20\,n\!+\!2\!-\!16\,nk\right )\left (n\!-
\!k\right )\left (a\!-\!k\!-\!1\right )\left (2\,a\!+\!k\right )}{n\left (2\,a\!+\!1\!+\!3\,
n\!-\!2\,k\right )\left (2\,a\!-\!1\!+\!3\,n\!-\!2\,k\right )\left (2\,a\!-\!3\!+\!3\,n\!-\!2\,k
\right )\left (n\!-\!1\right )}}
$
}
\\[3mm]
\hline \\[-3mm]
(1.7) & $1$ &
\multicolumn{2}{l}{
$\displaystyle
{\frac {\left (2\,a-1+4\,n+2\,d\right )\left (a-d+k\right )\left (2\,a
+k\right )\left (n-k\right )\left (2\,b-1-k\right )\left (2\,b+k
\right )}{n\left (2\,a\!+\!3\,k\right )\left (2\,d\!+\!2\,n\!-\!k\right )\left (2
\,d\!+\!2\,n\!-\!1\!-\!k\right )\left (2\,a\!+\!2\,b\!-\!1\!+\!2\,n\right )\left (a\!-\!b\!+\!n\right 
)}}
$
}
\\[3mm]
\hline \\[-3mm]
(1.8) & $2$ &
\multicolumn{2}{l}{
$\displaystyle
{8\frac {\left (2\,d\!-\!1\!-\!k\right )\left (b\!+\!k\right )\left (n\!-\!k\right 
)\left (a\!+\!k\right )\left (2\,a\!-\!2\,b\!+\!2\,k\!+\!1\right )\left (a\!+\!n\!+\!d\!-\!1
\right )}{n\left (2\,a\!+\!3\,k\right )\left (-2\!+\!2\,d\!+\!n\!-\!2\,k\right )\left 
(2\,b\!-\!1\!+\!n\right )\left (2\,a\!-\!2\,b\!+\!n\right )\left (2\,a\!+\!2\,d\!+\!n\!+\!k\!-\!1
\right )}}
$
}
\\[3mm]
\hline \\[-3mm]
(3.7) & $1$ &
\multicolumn{2}{l}{
$\displaystyle
{4\frac {\left (n-k\right )\left (6\,n+2-3\,k\right )\left (7\,n-
1-3\,k\right )}{\left (3\,n+1\right )\left (1+2\,n\right )n}}
$
}
\\[3mm]
\hline \\[-3mm]
(5.21) & $3$ &
\multicolumn{2}{l}{
$\displaystyle
{2\frac {\left (3\,a+1+k\right )\left (6\,a+2\,k+1\right )\left (n-k
\right )}{n\left (6\,a+n\right )\left (-n+6\,a+3+3\,k\right )}}
$
}
\\[3mm]
\hline \\[-3mm]
(5.22) & $1$ &
\multicolumn{2}{l}{
$\displaystyle
{\frac {\left (5+6\,k\right )\left (1+2\,k\right )\left (n-k\right )}{
\left (24\,n+4\right )\left (6\,n-1\right )n}}
$
}
\\[3mm]
\hline \\[-3mm]
(5.23) & $1$ &
\multicolumn{2}{l}{
$\displaystyle
{4\frac {\left (21\,n-7-9\,k\right )\left (6\,n+1-3\,k\right )\left (
n-k\right )}{\left (6\,n+1\right )\left (3\,n-1\right )n}}
$
}
\\[3mm]
\hline \\[-3mm]
(5.24) & $2$ &
\multicolumn{2}{l}{
$\displaystyle
{4\frac {n-k}{1+3\,n}}
$
}
\\[3mm]
\hline \\[-3mm]
(5.25) & $2$ &
\multicolumn{2}{l}{
$\displaystyle
{\frac {\left (4\,n-4\,k\right )\left (1+2\,k\right )\left (2+3\,k
\right )}{n\left (3\,n-1\right )\left (1+3\,n\right )}}
$
}
\\[3mm]
\hline \\[-3mm]
(5.27) & $1$ &
\multicolumn{2}{l}{
$\displaystyle
{81\frac {\left (n-1-2\,k\right )\left (n-2\,k\right )\left (-1+3
\,n-3\,k\right )}{n\left (3\,n-1\right )\left (-1+9\,n-21\,k\right )}}
$
}
\end{tabular}
\end{table}
}\noindent
with respect to the variable $k$. (In most cases $l=1$, so that
Gosper's original algorithm is applied.)
If successful, this generates $G(n,k)$ with
\be
a_k=F(n,k)-F(n-m,k)=G(n,k)-G(n,k-l)
\;,
\label{eq:WZ2extended}
\ee
and summing over all $k$ leads to
\[
s_n-s_{n-m}=\sum_{k\in\ZZ} \Big( F(n,k)-F(n-m,k)\Big)
=\sum_{k\in\ZZ} \Big( G(n,k)-G(n,k-l)\Big)=0
\]
since the right hand side is telescoping. Therefore $s_n$ is constant
for constant $n$ mod $m$, and these constants can be calculated using suitable
initial values. This can be accomplished if the series considered
is terminating. Note, that again, the function
\be
R(n,k)=\frac{G(n,k)}{F(n,k)}
\label{eq:rationalcerificateextended}
\ee
acts as a rational certificate function. Once the rational certificate
is known, it is a matter of pure rational arithmetic to decide the validity of
(\ref{eq:WZ1extended}) since the only thing that one has to show is
(\ref{eq:WZ2extended}) which after division by $F(n,k)$
is equivalent to the purely rational identity
\[
1-R(n,k)+R(n,k-l)\frac{F(n,k-l)}{F(n,k)}-\frac{F(n-m,k)}{F(n,k)}=0
\;.
\]
As an example, we prove (\ref{eq:Gessel-Stanton (1.1)}):
In the given case, we set $m:=3$, $l:=1$,
further
\[
F(n,k):=\frac{(-n)_k\,(n+3a)_k\,(a)_k}{k!\,(3a/2)_k\,((3a+1)/2)_k}
\frac{(n/3)!\,(3a+1)_n}{n!\,(a+1)_{n/3}}
\left(\frac{3}{4}\right)^k
\;,
\]
and notice that
\[
\frac{F(n,k)}{F(n,k-1)}
\quad\mbox{and}\quad
\frac{F(n,k)}{F(n-3,k)}
\]
are (complicated) rational functions 
(Algorithms~\ref{algo:1} and \ref{algo:3}). 
An application of Gosper's
algorithm is successful, and leads to the rational certificate
\[
R(n,k)=
{3\frac {\left (a+k\right )\left (n-k\right )\left (3\,a+2\,n-3
\right )}{\left (n+3\,a+k-2\right )\left (n+3\,a+k-1\right )n}}
\;.
\]
Therefore
\[
\sum\limits_{k\in\ZZ} F(n,k)=\sum\limits_{k=0}^n F(n,k)=
\mbox{constant}\quad\quad{(n \mbox{ mod }} 3 \mbox{ constant})
\;,
\]
and statement (\ref{eq:Gessel-Stanton (1.1)})
follows using three trivial initial values.

Table \ref{table:Gessel-Stanton} lists the
hypergeometric identities of the Gessel-Stanton paper
(note the misprint in Equation (1.4)), and
Table \ref{table:WZextended} contains their rational certificates
(\ref{eq:rationalcerificateextended}), calculated
by our implementations,
together with the certificates of Bailey's list (Table \ref{table:Bailey})
to which the WZ method did not apply. 

Note that in all cases considered, $l=1$, so that the original Gosper
algorithm is applied.

Note, moreover, that Gessel and Stanton were not able to present proofs for
their statements (6.2), (6.3), (6.5), and (6.6):
Table \ref{table:Gessel and Stanton's open problems} contains proofs.

{\small
\begin{table}[hbtp]
\caption{Gessel and Stanton's open problems}
\label{table:Gessel and Stanton's open problems}
\vspace*{3mm}
\begin{tabular}{lrp{10cm}}
Eq. & \multicolumn{2}{l}{
Identity}
\\[2mm] \hline \\[-2mm]
(6.2)
&
\multicolumn{2}{l}{
$\displaystyle
_7 F_6\left.
\!\!
\left(
\!\!\!\!
\begin{array}{c}
\multicolumn{1}{c}{\begin{array}{c}
a+1/2\;, a\;, b\;, 1-b\;, -n\;,(2a+1)/3+n\;, a/2+1
\end{array}}\\[1mm]
\multicolumn{1}{c}{\begin{array}{c}
1/2\;, (2a-b+3)/3\;, (2a+b+2)/3\;, -3n\;, 2a+1+3n\;, a/2
            \end{array}}\end{array}
\!\!\!\!
\right| 1\right)
=$
}
\\[4mm]
&
\multicolumn{2}{l}{
$
\hspace*{1.5cm}
\displaystyle
\frac{((2a+2)/3)_n\,(2a/3+1)_n\,((1+b)/3)_n\,((2-b)/3)_n}
{((2a-b)/3+1)_n\,((2a+b+2)/3)_n\,(2/3)_n\,(1/3)_n}
$
}
\\[5mm]
(6.3)
&
\multicolumn{2}{l}{
$\displaystyle
_5 F_4\left.
\!\!
\left(
\!\!\!\!
\begin{array}{c}
\multicolumn{1}{c}{\begin{array}{c}
a+1/2\;, a\;, -n\;,(2a+1)/3+n\;, a/2+1
\end{array}}\\[1mm]
\multicolumn{1}{c}{\begin{array}{c}
1/2\;, -3n\;, 2a+1+3n\;, a/2
            \end{array}}\end{array}
\!\!\!\!
\right| 9\right)
=
\frac{((2a+2)/3)_n\,(2a/3+1)_n}
{(2/3)_n\,(1/3)_n}
$
}
\\[5mm]
(6.5)
&
\multicolumn{2}{l}{
$\displaystyle
_2 F_1\left.
\!\!
\left(
\!\!\!\!
\begin{array}{c}
\multicolumn{1}{c}{\begin{array}{c}
-n\;,-n+1/4
\end{array}}\\[1mm]
\multicolumn{1}{c}{\begin{array}{c}
2n+5/4
            \end{array}}\end{array}
\!\!\!\!
\right| \frac{1}{9}\right)
=
\frac{(5/4)_{2n}}{(2/3)_n\,(13/12)_n}\left(\frac{2^6}{3^5}\right)^n
$
}
\\[5mm]
(6.6)
&
\multicolumn{2}{l}{
$\displaystyle
_2 F_1\left.
\!\!
\left(
\!\!\!\!
\begin{array}{c}
\multicolumn{1}{c}{\begin{array}{c}
-n\;,-n+1/4
\end{array}}\\[1mm]
\multicolumn{1}{c}{\begin{array}{c}
2n+9/4
            \end{array}}\end{array}
\!\!\!\!
\right| \frac{1}{9}\right)
=
\frac{(9/4)_{2n}}{(4/3)_n\,(17/12)_n}\left(\frac{2^6}{3^5}\right)^n
$
}
\\[2mm] \hline \\[-2mm]
\multicolumn{3}{c}{Rational certificates}\\[2mm]
\hline \\[-2mm]
Eq.&$m$&$R(n,k)$
\\[2mm] \hline \\[-2mm]
(6.2) & $1$ &
$\displaystyle
{6\frac {\left (a-1+3\,n\right )\left (a+k\right )\left (2\,a+2\,k+
1\right )\left (n-k\right )\left (b-1-k\right )\left (b+k\right )}{
\left (a+2\,k\right )\left (3\,n-k\right )\left (3\,n-1-k\right )
\left (3\,n-2-k\right )\left (2\,a-b+3\,n\right )\left (2\,a+b-1+3\,n
\right )}}
$
\\[2mm] \hline \\[-2mm]
(6.3) & $1$ &
$\displaystyle
-{\frac {\left (6\,a-6+18\,n\right )\left (n-k\right )\left (2\,a+2\,k
+1\right )\left (a+k\right )}{\left (a+2\,k\right )\left (3\,n-k
\right )\left (3\,n-1-k\right )\left (3\,n-2-k\right )}}
$
\\[2mm] \hline \\[-2mm]
(6.5) & $1$ &
$\displaystyle
-{\frac {\left (52\,{n}^{2}-13\,n-21-56\,k+16\,nk-32\,{k}^{2}\right )
\left (n-k\right )\left (4\,n-1-4\,k\right )}{\left (108\,n-27\right )
\left (3\,n-1\right )\left (1+12\,n\right )n}}
$
\\[2mm] \hline \\[-2mm]
(6.6) & $1$ &
$\displaystyle
-{\frac {\left (52\,{n}^{2}+39\,n-55-84\,k+16\,nk-32\,{k}^{2}\right )
\left (4\,n-1-4\,k\right )\left (n-k\right )}{\left (108\,n-27\right )
\left (1+3\,n\right )\left (5+12\,n\right )n}}
$

\end{tabular}
\end{table}
}\noindent

Similarly as the original WZ approach, our method is not capable, however, 
to prove Gessel-Stanton's (6.1), a non-terminating version of (6.2). 
Also, Gessel-Stanton's result (1.9)
\[
_3 F_2\left.
\!\!
\left(
\!\!\!\!
\begin{array}{c}
\multicolumn{1}{c}{\begin{array}{c}
-sb+s+1\;, b-1\;,-n
\end{array}}\\[1mm]
\multicolumn{1}{c}{\begin{array}{c}
b+1\;, s(-n-b)-n
            \end{array}}\end{array}
\!\!\!\!
\right| 1\right)
=
\frac{(1+s+sn)_{n}\,b\,(n+1)}{(1+s(b+n))_n\,(b+n)}
\]
is over the capabililites of our method since in this case the summand
is an $(m,l)$-fold hypergeometric term only for fixed (rational), 
but not for arbitrary $s$.

Note, that our method not only unifies the proof of hypergeometric identities
in a stronger fashion than the original WZ approach but moreover
our \Reduce\ and \Maple\ implementations do all the necessary computations
completely automatically.
We present some of the input and output in the appendix.

Finally, we give examples of an application for which $l\neq 1$.
To prove the identity $(n\in\N)$
\be
-\sum_{k=0}^n (-2)^n\,\ueber{n}{k} \cdot \ueber{k/2}{n}=1
\;,
\label{eq:k/2}
\ee
we apply our extended WZ method with $l=2$, and $m=1$,
and get the rational certificate
\[
R(n,k)=
{\frac {\left (-k+n-1\right )\left (-k+n\right )}{\left (n-1\right )
\left (-k+2\,n-2\right )}}
\;,
\]
which proves (\ref{eq:k/2}). Similarly one proves the statement $(n\in\N_0)$
\[
\sum_{k=0}^n (-1)^k\,(-2)^n\,\ueber{n}{k} \cdot \ueber{k/2}{n}=1
\;.
\]

\section{The Zeilberger Algorithm}
\label{sec:The Zeilberger Algorithm}

In this section, we recall the celebrated Zeilberger algorithm 
(Zeilberger, 1990--1991), see also Graham, Knuth and Patashnik (1994)
with which one can not only verify hypergeometric identities
but moreover definite sums can be calculated if they represent
hypergeometic terms.

Zeilberger's algorithm deals with the question to determine a
{\sl holonomic recurrence equation}
\be
\sum_{j=0}^{J} P_j(n)\,\Sigma(n-j)=0
\label{eq:zeilbergerrecurrence}
\ee
with polynomials $P_j$ in $n$, for sums
\be
\Sigma(n):=\sum_{k\in\ZZ} F(n,k)
\label{eq:zeilbergersumme}
\ee
for which $F(n,k)$ is a hypergeometric term with respect to both $n$ and $k$.

Zeilberger's idea is to apply Gosper's algorithm in the following
non-obvious way: Set 
\[\
a_k:=F(n,k)+\sum_{j=1}^J \sigma_j(n)\,F(n-j,k)
\]
with yet undetermined variables $\sigma_j$ depending on $n$, but not
depending on $k$. Then
\[
\frac{a_k}{a_{k-1}}=
\frac{F(n,k)+\sum\limits_{j=1}^J \sigma_j(n)\,F(n-j,k)}
{F(n,k-1)+\sum\limits_{j=1}^J \sigma_j(n)\,F(n-j,k-1)}
=
\frac{F(n,k)}{F(n,k-1)}\cdot
\frac{1+\sum\limits_{j=1}^J \sigma_j(n)\,\frac{F(n-j,k)}{F(n,k)}}
{1+\sum\limits_{j=1}^J \sigma_j(n)\,\frac{F(n-j,k-1)}{F(n,k-1)}}
\]
turns out to be rational with respect to $k$, so the Gosper algorithm
may be applied.

If an application of Gosper's algorithm is successful it provides us with
$s_k$ depending on $n$, and a set of rational functions $\sigma_j(n)$ 
(the coefficients of $f_k$ are determined together with the unknowns
$\sigma_j$) such that
\[
s_k-s_{k-1}=a_k=F(n,k)+\sum_{j=1}^J \sigma_j(n)\,F(n-j,k)
\;,
\]
so that by summation
\beao
\sum_{k\in\ZZ}a_k
&=&
\sum_{k\in\ZZ}\lk F(n,k)+\sum_{j=1}^J \sigma_j(n)\,F(n-j,k)\rk
\\&=&
\Sigma(n)+\sum_{j=1}^J \sigma_j(n)\,\Sigma(n-j)
=
\sum_{k\in\ZZ}\Big( s_k-s_{k-1}\Big)=0
\eeao
since the right hand side is a telescoping sum. After multiplication with the
common denominator this establishes the recurrence
equation (\ref{eq:zeilbergerrecurrence}) searched for.

Koornwinder (1993) gives a rigorous description of Zeilberger's
algorithm in the (most common) case that the summation bounds are natural:
$a_{-1}=a_{n+1}=0$, i.\ e.\ the summation is for $k=0\ldots n$.

Like for the Wilf-Zeilberger method, the Zeilberger algorithm is
accompanied by a rational certification mechanism.

Note that Zeilberger's algorithm can be applied to ratios of products 
of rational functions, exponentials, factorials, $\Gamma$ function terms, binomial
coefficients, and Pochhammer symbols that are integer-linear in their
arguments with respect to both $n$ and $k$. 

In the next section we will present a modified version of
Zeilberger's algorithm that is applicable if the arguments of
such expressions are rational-linear with respect to $n$ and $k$.

The application of Zeilberger's algorithm has the advantage 
over the WZ method that the right hand side
of the hypergeometric identity does not have to be known in advance,
but is {\sl generated} by the algorithm (not to speak of the possibility
to verify identities of other type).
Therefore, Zeilberger's algorithm can be used to calculate definite sums
rather than only verifying them.
All identities mentioned in this article which could be verified
with the WZ method, can be generated by Zeilberger's algorithm.

Implementations of the Zeilberger
algorithm were given by Zeilberger (1990) and Koornwinder (1993) in
\Maple, and by Paule and Schorn (1994) in \Mathematica.
On the lines of Koornwinder (1993), we implemented the Zeilberger algorithm
in \Reduce\ (Koepf, 1994) and \Maple,
examples of which are given in the appendix.

Note the following side conditions of the previous implementations:
\begin{itemize}
\item
Zeilberger: Here one must write the input into a file rather than on
the command line. Supports only integer-linear input of a special
form.
\item
Koornwinder: Supports only integer-linear input in hypergeometric notation.
\item
Paule-Schorn: Supports only ratios of rational functions, 
products of exponentials, factorials, and binomial coefficients.
\end{itemize}
Our implementations support the input in 
factorial-binomial-Gamma-Pochhammer as well as hypergeometric 
notation, and use Algorithm~\ref{algo:3} for rationality decisions, and 
are therefore not bound to integer-linear input.


\section{An Extended Version of Zeilberger's Algorithm}
\label{sec:An Extended Version of Zeilberger's Algorithm}

%

Our extended version of the Zeilberger algorithm deals with the
question to determine a holonomic recurrence equation
(\ref{eq:zeilbergerrecurrence}) for sums (\ref{eq:zeilbergersumme})
for which $F(n,k)$ is an $(m,l)$-fold hypergeometric term with respect to
$(n,k)$, see \S~\ref{sec:Extension of the WZ method}.

In particular, this applies to all cases when the input function $F(n,k)$
is given as a ratio of products of
rational functions, exponentials, factorials, $\Gamma$ function terms, binomial
coefficients, and Pochhammer symbols that are rational-linear in their
arguments with respect to both $n$, and $k$.

First of all we mention that Zeilberger's algorithm may be applicable
even though this is safely the case only if the arguments are integer-linear.
An example of that type is the function
\[
\Sigma(n):=\;
_2 F_1\left.
\!\!
\left(
\!\!\!\!
\begin{array}{c}
\multicolumn{1}{c}{\begin{array}{cc} -n/2\;, & -n/2+1/2 \end{array}}\\[1mm]
\multicolumn{1}{c}{ b+1/2}
            \end{array}
\!\!\!\!
\right| 1\right)
=
\sum_{k=0}^\infty \frac{(-n/2)_k\,(-n/2+1/2)_k}{k!\,(b+1/2)_k}
\;,
\]
for which
an application of Zeilberger's algorithm yields the recurrence equation
\[
(2b + n-1)\,\Sigma(n) - 2(b + n-1)\,\Sigma(n-1)=0
\;,
\]
and therefore the explicit representation
\[
\Sigma(n)=\frac{2^n\,(b)_n}{(2b)_n}
\;.
\]
Zeilberger's algorithm applies since $F(n,k)/F(n-1,k)$ and
$F(n,k)/F(n,k-1)$ are rational even
though the representing expression for $F(n,k)$ is not integer-linear in
its arguments.

On the other hand, not for every $F(n,k)$ given with
rational-linear $\Gamma$-arguments, the Zeilberger algorithm is
applicable. An example for this situation is the left hand side
of the Watson theorem with respect to variable $a$ 
(see Table~\ref{table:Bailey}).

We present now an algorithm which can be applied for arbitrary 
rational-linear input.
\begin{algorithm}
\label{algo:4}
{\rm
({\tt extended\verb+_+sumrecursion})\\[3mm]
The following steps perform an algorithm to
determine a holonomic recurrence equation
(\ref{eq:zeilbergerrecurrence}) for sums (\ref{eq:zeilbergersumme}).
\begin{enumerate}
\item
Input: $F(n,k)$, given as a ratio of products of
rational functions, exponentials, factorials, $\Gamma$ function terms, binomial
coefficients, and Pochhammer symbols
with rational-linear arguments in $n$ and $k$. 
\item
Build the list of all arguments. They are of the form 
$p_j/q_j\,n+s_j/t_j\,k+\al_j$ with integer $p_j, q_j, s_j, t_j$,
$p_j/q_j$ and $s_j/t_j$ in lowest terms, $q_j$ and $t_j$ positive.
\item
Calculate $m:=\lcm\{q_j\}$ and $l:=\lcm\{t_j\}$.
\item
Define $\tilde F(n,k):=F(mn,kl)$. 
Then $\tilde F(n,k)$ is integer-linear
in the arguments.
\item
Apply Zeilberger's algorithm to $\tilde F(n,k)$.
Get the recurrence equation
\[
\sum_{j=0}^{J} P_j(n)\,\tilde\Sigma(n-j)=0
\]
with polynomials $P_j$ in $n$, for the sum
\[
\tilde\Sigma(n):=\sum_{k\in\ZZ} \tilde F(n,k)
\;.
\]
\item
The output is the recurrence equation
\[
\sum_{j=0}^{J} P_j(n/m)\,
\Sigma(n-jm)=0
\]
for the sum
\[
\Sigma(n):=\sum_{k\in\ZZ} F(n,k)
\;.
\]
\end{enumerate}
}
\end{algorithm}
\begin{proof}
Obviously our construction provides us with $\tilde F(n,k)$ 
integer-linear in the arguments involved. Therefore
Zeilberger's algorithm can be applied, and the result follows.
\end{proof}
\\[-3mm]
As a first example, we apply our algorithm to the Watson function
\[
\Sigma(n)
=\;
_3 F_2\left.
\!\!
\left(
\!\!\!\!
\begin{array}{c}
\multicolumn{1}{c}{\begin{array}{c}
-n\;, b\;, c
\end{array}}\\[1mm]
\multicolumn{1}{c}{\begin{array}{c}
(-n+b+1)/2 \;, 2c
            \end{array}}\end{array}
\!\!\!\!
\right| 1\right)
\]
with respect to the variable $n$
to which Zeilberger's algorithm does not apply. In this case, the algorithm
determines $m=2$ and $l=1$, and leads to the two-fold recurrence equation
\[
(b - 2c - n + 1)\,(n - 1)\,\Sigma(n - 2) - (b - n + 1)\,(2c + n - 1)\,\Sigma(n)
=0
\]
from which the explicit right hand representation listed in 
Table~\ref{table:Bailey}
can be deduced for integer $n$ since for positive values of $n$ the 
Watson sum is finite, and therefore
\[
\Sigma(0)=
1
\;,
\]
and
\[
\Sigma(1)=
1+\frac{-1\,b\,c}{1\,(b/2)\,(2c)}=0
\;.
\]
It turns out that our method is applicable to all identities considered
in this paper to which Zeilberger's original approach does not apply.

For example, we consider the three major identities of
the paper of Gessel and Stanton (1982):
The evaluation of (1.7)
\beao
\Sigma(n)&:=&
_7 F_6\left.
\!\!
\left(
\!\!\!\!
\begin{array}{c}
\multicolumn{1}{c}{\begin{array}{c}
2a\;, 2b\;, 1-2b\;, 1-2a/3\;, a+d+n+1/2\;, a-d\;, -n
\end{array}}\\[1mm]
\multicolumn{1}{c}{\begin{array}{c}
a-b+1\;, a+b+1/2\;, 2a/3\;, -2d-2n\;, 2d+1\;, 1+2a+2n
            \end{array}}\end{array}
\!\!\!\!
\right| 1\right)
\\&=&
\frac{(a+1/2)_n\,(a+1)_n\,(b+d+1/2)_n\,(d-b+1)_n}
{(a+b+1/2)_n\,(a-b+1)_n\,(d+1/2)_n\,(d+1)_n}
\eeao
is found by a direct application of Zeilberger's algorithm with respect
to $n$, leading to the equivalent recurrence equation
\beao
0&=&
(2a + 2b + 2n - 1)(a - b + n)(2d + 2n - 1)(d + n)\Sigma(n)
\\&&
  + (2a + 2n - 1)(a + n)(2b + 2d + 2n - 1)(b - d - n)\Sigma(n - 1)
\;.
\eeao
On the other hand, the evaluation of (1.8)
\beao
\Sigma(n)&:=&
_7 F_6\left.
\!\!
\left(
\!\!\!\!
\begin{array}{c}
\multicolumn{1}{c}{\begin{array}{c}
a\;, b\;, a+1/2-b\;, 1+2a/3\;, 1-2d\;,2a+2d+n\;, -n
\end{array}}\\[1mm]
\multicolumn{1}{c}{\begin{array}{c}
2a-2b+1\;, 2b\;, 2a/3\;, a+d+1/2\;, 1-d-n/2\;, 1+a+n/2
            \end{array}}\end{array}
\!\!\!\!
\right| 1\right)
\\&=&
\funkdef{0}{n\;\mbox{odd}}{\displaystyle
\frac{(b+d)_{n/2}\,(d-b+a+1/2)_{n/2}\,n!\,(a+1)_{n/2}}
{(b+1/2)_{n/2}\,(a+d+1/2)_{n/2}\,(d)_{n/2}\,(n/2)!\,(a-b+1)_{n/2}}
}
\eeao
cannot be handled with respect to $n$ using Zeilberger's algorithm,
but the extended version leads to the equivalent 2-fold recurrence equation
\beao
0&=&
   (n - 1 + 2d + 2a)(2b - n - 2a)(n - 1 + 2b)(n - 2 + 2d)\Sigma(n)
\\&&
  +
(n - 1 + 2d - 2b + 2a)(n - 2 + 2d + 2b)(2a + n)(n - 1)\Sigma(n - 2)
\;.
\eeao
A direct application of Zeilberger's algorithm is possible,
however, with respect to the other variables involved 
(even with respect to $a$).

Gessel-Stanton's open problem (6.2)
\beao
\Sigma(n)&:=&
_7 F_6\left.
\!\!
\left(
\!\!\!\!
\begin{array}{c}
\multicolumn{1}{c}{\begin{array}{c}
a+1/2\;, a\;, b\;, 1-b\;, -n\;,(2a+1)/3+n\;, a/2+1
\end{array}}\\[1mm]
\multicolumn{1}{c}{\begin{array}{c}
1/2\;, (2a-b+3)/3\;, (2a+b+2)/3\;, -3n\;, 2a+1+3n\;, a/2
            \end{array}}\end{array}
\!\!\!\!
\right| 1\right)
\\&=&
\frac{((2a+2)/3)_n\,(2a/3+1)_n\,((1+b)/3)_n\,((2-b)/3)_n}
{((2a-b)/3+1)_n\,((2a+b+2)/3)_n\,(2/3)_n\,(1/3)_n}
\;,
\eeao
again, can be solved directly with Zeilberger's algorithm leading to
the recurrence equation
\beao
0&=&
(2a + b + 3n - 1)\,(2a - b + 3n)\,(3n - 1)\,(3n - 2)\,\Sigma(n)
\\&&
+(2a + 3n - 1)\,(2a + 3n)\,(b + 3n - 2)\,(b - 3n + 1)\,\Sigma(n - 1)
\;.
\eeao

Finally, as an example with $l\neq 1$, we consider (\ref{eq:k/2}), again.
Our algorithm generates $m=1$ and $l=2$, and the recurrence equations
\[
\Sigma(n) - \Sigma(n - 1)=0
\quad\quad\mbox{and}\quad\quad
2 \Sigma(n) + \Sigma(n - 1)=0
\]
for 
\[
\Sigma(n):=(-2)^n\,\ueber{n}{k} \cdot \ueber{k/2}{n}\;,
\quad\quad\mbox{and}\quad\quad
\Sigma(n):=\ueber{n}{k} \cdot \ueber{k/2}{n}\;,
\]
respectively.

\section{Deduction of hypergeometric identities}

Finally, we mention that with a good implementation of Zeilberger's
algorithm and our extension at hand, it is easy to {\sl discover}
new identities.
Just for fun, we realized the pattern in Andrews' statement 
(\ref{eq:Gessel-Stanton (1.1)}), and tried to generate similar ones:
It turns out that
\[
_3 F_2\left.
\!\!
\left(
\!\!\!\!
\begin{array}{c}
\multicolumn{1}{c}{\begin{array}{c}
-n\;, n+2a\;, a
\end{array}}\\[1mm]
\multicolumn{1}{c}{\begin{array}{c}
2a/2\;,(2a+1)/2
            \end{array}}\end{array}
\!\!\!\!
\right| \frac{2}{4}\right)
=
\funkdef{0}{n {\mbox{ odd }}}{\displaystyle
\frac{(-1)^{n/2}\,(1/2)_{n/2}}
{(1/2+a)_{n/2}}
}
\]
and
\[
_3 F_2\left.
\!\!
\left(
\!\!\!\!
\begin{array}{c}
\multicolumn{1}{c}{\begin{array}{c}
-n\;, n+4a\;, a
\end{array}}\\[1mm]
\multicolumn{1}{c}{\begin{array}{c}
4a/2\;,(4a+1)/2
            \end{array}}\end{array}
\!\!\!\!
\right| \frac{4}{4}\right)
=
\funkdef{0}{n {\mbox{ odd }}}{\displaystyle
\frac{(1/2)_{n/2}}
{(1/2+2a)_{n/2}}
}
\;.
\]
Another example of a more deductive strategy is: Applying Zeilberger's
algorithm to the general $_2 F_1$ polynomial
\[
\Sigma(n):=\;
_2 F_1\left.
\!\!
\left(
\!\!\!\!
\begin{array}{c}
\multicolumn{1}{c}{\begin{array}{c}
a\;, -n
\end{array}}\\[1mm]
\multicolumn{1}{c}{\begin{array}{c}
b
            \end{array}}\end{array}
\!\!\!\!
\right| x\right)
\;,
\]
e.\ g., leads to the recurrence equation
\[
(b - 1 + n)\,\Sigma(n)+(- 2 n + x n + x a - b + 2 - x)\,\Sigma(n-1)-
(x - 1)\, (n - 1)\,\Sigma(n-2) =0
\;.
\]
It is therefore hypergeometric only if the coefficient of $\Sigma(n-2)$
is identical zero, 
i.\ e.\ if $x=1$. This gives Vandermonde's identity. However,
the coefficient of $\Sigma(n-1)$ can be made zero (equating
coefficients), if we choose $x=2$,
and $b=2a$, in which situation we get
\[
(n + 2 a - 1)\,\Sigma(n)- (n - 1)\, \Sigma(n-2)=0
\;.
\]
Therefore we have deduced the identity
\[
_2 F_1\left.
\!\!
\left(
\!\!\!\!
\begin{array}{c}
\multicolumn{1}{c}{\begin{array}{c}
a\;, -n
\end{array}}\\[1mm]
\multicolumn{1}{c}{\begin{array}{c}
2a
            \end{array}}\end{array}
\!\!\!\!
\right| 2\right)
=
\funkdef{0}{n {\mbox{ odd }}}{\displaystyle
\frac{(1/2)_{n/2}}
{(1/2+a)_{n/2}}
}
\;.
\]
We see that this method, to some extent,
can be a substitute for the ingenuity of people like
Dougall, Bailey, Andrews, Gessel or Stanton to find hypergeometric
sums which can be represented by single hypergeometric terms.

We finally give a strange example to demonstrate
that our method can be of great help to find new identities.

We try to find all hypergeometric functions of the form
\[
\Sigma(n):=\;
_2 F_1\left.
\!\!
\left(
\!\!\!\!
\begin{array}{c}
\multicolumn{1}{c}{\begin{array}{c}
a\;, -n
\end{array}}\\[1mm]
\multicolumn{1}{c}{\begin{array}{c}
n+b
            \end{array}}\end{array}
\!\!\!\!
\right| x\right)
\]
for which $a, b$ and $x$ are constants with respect to $n$, and for which
a recurrence equation with only two terms $\Sigma(n-j)$ is valid.

The recurrence equation for $\Sigma(n)$ turns out to be
\begin{eqnarray*}
0&=&
-\left (x-1\right )^{2}\left (n-1\right )\left (n-1+b\right )\left (n-
2+b\right )\left (xn+n-xa-x+bx\right )\Sigma(n-2)
\\&&
+\left (n-1+b \right )
P(n,a,b,x)\,\Sigma(n-1)
\\&&
+ x\left (2\,n+b-1\right )\left (2\,n+b-2
\right )\left (n-a-1+b\right )\left (xn+n-xa-2\,x-1+bx\right )
\Sigma(n)
\;,
\end{eqnarray*}
where $P(n,a,b,x)$ denotes a very complicated
polynomial of degree $2$ in $n$ that does not have a rational factorization.
To receive a recurrence equation for which only two terms $\Sigma(n-j)$
different from zero occur, we may set the coefficient lists with
respect to $n$ of any of the factors occurring zero, and try to solve
for $a,b$ and $x$.
Note that since the resulting equations systems are polynomial systems, by
Gr\"obner bases methods these can be solved algorithmically. 

In our case, we receive either $x=1$, i.\ e.\ the recurrence equation 
\begin{eqnarray*}
0&=&
-\left (n-1+b\right )\left (2\,n-a+b-1\right )\left (b-a+2\,n-2\right 
)\Sigma(n-1)
\\&&
+\left (2\,n+b-1\right )\left (2\,n+b-2
\right )\left (n-a-1+b\right )\Sigma(n)
\;,
\end{eqnarray*}
or we are led to the Kummer identity, i.\ e.\ to the values
$b=a+1$ and $x=-1$ with the recurrence equation
\[
-2\,\left (n+a\right )\Sigma(n-1)+\left (2\,n+a\right )\Sigma(n)
=0
\;.
\]
The only exception occurs when we set the coefficient list with
respect to $n$ of the factor $P(n,a,b,x)$ zero, leading to the Kummer
case again, and to the second solution set
\[
\{a=1/2,b=3/2,x^2-6x+1=0\}
\;.
\]
For $x=3\pm 2\sqrt 2$, we have the recurrence equation
\[
-4\,\left (2\,n-1\right )\left (2\,n+1\right )\Sigma(n-2)+
\left (4\,n-1\right )\left (4\,n+1\right )\Sigma(n)
=0
\]
leading to the closed form representations
\[
_2 F_1\left.
\!\!
\left(
\!\!\!\!
\begin{array}{c}
\multicolumn{1}{c}{\begin{array}{c}
1/2\;, -n
\end{array}}\\[1mm]
\multicolumn{1}{c}{\begin{array}{c}
n+3/2
            \end{array}}\end{array}
\!\!\!\!
\right| 3+2\sqrt 2\right)
=
\funkdef{\displaystyle
\frac{2\,(5/4)_{(n-1)/2}\,(7/4)_{(n-1)/2}}
{5\,(11/8)_{(n-1)/2}\,(13/8)_{(n-1)/2}}\, (1\!-\!\sqrt 2)
}{n {\mbox{ odd }}}{\displaystyle
\frac{(3/4)_{n/2}\,(5/4)_{n/2}}
{(7/8)_{n/2}\,(9/8)_{n/2}}
}
\]
and
\[
_2 F_1\left.
\!\!
\left(
\!\!\!\!
\begin{array}{c}
\multicolumn{1}{c}{\begin{array}{c}
1/2\;, -n
\end{array}}\\[1mm]
\multicolumn{1}{c}{\begin{array}{c}
n+3/2
            \end{array}}\end{array}
\!\!\!\!
\right| 3-2\sqrt 2\right)
=
\funkdef{\displaystyle
\frac{2\,(5/4)_{(n-1)/2}\,(7/4)_{(n-1)/2}}
{5\,(11/8)_{(n-1)/2}\,(13/8)_{(n-1)/2}}\, (1\!+\!\sqrt 2)
}{n {\mbox{ odd }}}{\displaystyle
\frac{(3/4)_{n/2}\,(5/4)_{n/2}}
{(7/8)_{n/2}\,(9/8)_{n/2}}
}
\;,
\]
in particular, for even $n$, the values at $x=3+2\sqrt 2$ and $x=3-2\sqrt 2$
are rational and equal:
\[
_2 F_1\left.
\!\!
\left(
\!\!\!\!
\begin{array}{c}
\multicolumn{1}{c}{\begin{array}{c}
1/2\;, -2n
\end{array}}\\[1mm]
\multicolumn{1}{c}{\begin{array}{c}
2n+3/2
            \end{array}}\end{array}
\!\!\!\!
\right| 3\pm 2\sqrt 2\right)
=
\sum_{k=0}^{2n} (-1)^k\,\frac{\ueber{2n}{k}\,\ueber{2n+k+1}{k}}
{\ueber{4n+2k+2}{2k}} \lk 3\pm 2\sqrt 2\rk^k
=
\frac{(3/4)_{n}\,(5/4)_{n}}
{(7/8)_{n}\,(9/8)_{n}}
\;.
\]
We will discuss the given method in greater detail in a forthcoming paper.
 
\section*{Acknowledgement}
I like to thank Gregor St\"olting for his work on the implementations,
Tom Koornwinder for his wonderful implementation {\tt zeilb} which was
the starting point of our implementations, and Peter Deuflhard who
initiated my studies on the given topic.

\section*{Appendix}

In this appendix, we give a short description of a \Maple\ implementation,
which I implemented together with Gregor St\"olting on the lines of
Koornwinder (1993),
incorporating Gosper's and Zeilberger's algorithms and the extensions
of this article, and present some of its results.
Our \Reduce\ implementation is described elsewhere
(Koepf, 1994).

After loading our package, one can use the following \Maple\ functions:
\begin{description}
\vspace*{-3mm}
\item
{\tt gosper(f,k)} determines a closed
form antidifference. If it does not return a closed form solution, then
a closed form solution does not exist.
\item
{\tt gosper(f,k,m,n)} determines
\[
\sum_{k=m}^n a_k
\]
using Gosper's algorithm. This is only successful if Gosper's algorithm applies.
\item
{\tt extended\verb+_+gosper(f,k,m)} determines an  $m$-fold hypergeometric
antidifference. If it does not return a solution, then
such a solution does not exist.
\item
{\tt sumrecursion(f,k,n)} determines a holonomic recurrence equation
for \\
${\tt summ(n)}
=\sum\limits_{k=-\infty}^\infty f(n,k)$ with respect to $n$ if $f(n,k)$
is hypergeometric with respect to both $n$ and $k$.
The resulting expression equals zero.
\item
{\tt sumrecursion(f,k,n,j)} 
searches only for a holonomic recurrence equation of order $j$.
\item
{\tt extended\verb+_+sumrecursion(f,k,n,m,l)} determines a holonomic recurrence 
equation for \\
${\tt summ(n)}
=\sum\limits_{k=-\infty}^\infty f(n,k)$ with respect to $n$ if $f(n,k)$
is an $(m,l)$-fold hypergeometric term with respect to $(n,k)$.
\item
{\tt hyperrecursion(upper,lower,x,n)} determines a holonomic recurrence
equation with respect to $n$ for
$_{p}F_{q}\left.\left(\begin{array}{cccc}
a_{1},&a_{2},&\cdots,&a_{p}\\
b_{1},&b_{2},&\cdots,&b_{q}\\
            \end{array}\right| x\right)
$, where {\tt upper} $=\{a_{1}, a_{2}, \ldots, a_{p}\}$
is the list of upper parameters, and
{\tt lower} $=\{b_{1}, b_{2}, \ldots, b_{q}\}$
is the list of lower parameters depending on $n$.
\item
{\tt hyperrecursion(upper,lower,x,n,j)} 
searches only for a holonomic recurrence equation of order $j$.
\item
{\tt hyperterm(upper,lower,x,k)} yields the hypergeometric term
\[
\frac
{(a_{1})_{k}\cdot(a_{2})_{k}\cdots(a_{p})_{k}}
{(b_{1})_{k}\cdot(b_{2})_{k}\cdots(b_{q})_{k}\,k!}x^{k}
\]
with upper parameters {\tt upper} $=\{a_{1}, a_{2}, \ldots, a_{p}\}$,
and lower parameters {\tt lower} $=$\linebreak
$\{b_{1}, b_{2}, \ldots, b_{q}\}$
\item
{\tt simplify\verb+_+gamma(f)} simplifies an expression {\tt f} involving
only rational functions, exponentials and
$\Gamma$ function terms according to a recursive
application of the simplification rule\linebreak
$\Gamma\:(a+1)=a\,\Gamma\:(a)$
to the expression tree, see Algorithm~\ref{algo:3}.
\item
{\tt simplify\verb+_+combinatorial(f)} simplifies an expression {\tt f}
involving exponentials, factorials, $\Gamma$ function terms,
binomial coefficients, and Pochhammer symbols by converting
factorials, binomial coefficients, and Pochhammer symbols into
$\Gamma$ function terms, and applying {\tt simplify\verb+_+gamma} 
and {\tt simplify\verb+_+power} to its result. If the output is not rational,
it is given in terms of $\Gamma$ functions, see Algorithm~\ref{algo:3}.
\vspace*{-3mm}
\end{description}
The \Maple\ function

\begin{verbatim}
WZ:=proc(summand,k,n,m)
local tmp,gos;
   tmp:=summand-subs(n=n-m,summand);
   gos:=extended_gosper(tmp,k,m);
   RETURN(simplify_combinatorial(gos/summand))
end:
\end{verbatim}
therefore, calculates the $(m,l)$-fold rational certificate 
(\ref{eq:rationalcerificateextended}) of $F(n,k)$.

Here are some results of the implementation:

{\small
\begin{verbatim}
    |\^/|     Maple V Release 3 (FU-Berlin)
._|\|   |/|_. Copyright (c) 1981-1994 by Waterloo Maple Software and the
 \  MAPLE  /  University of Waterloo. All rights reserved. Maple and Maple V
 <____ ____>  are registered trademarks of Waterloo Maple Software.
      |       Type ? for help.
> read summation;

> # see (SIAM Review, 1994, Problem 94-2)

> gosper((-1)^(k+1)*(4*k+1)*(2*k)!/(k!*4^k*(2*k-1)*(k+1)!),k);
                                          (k + 1)
                               (2 k)! (-1)
                               ------------------
                                           k
                                 (k + 1)! 4  k!
> # Dougall

> WZ(hyperterm({a,1+a/2,b,c,d,1+2*a-b-c-d+n,-n},
  {a/2,1+a-b,1+a-c,1+a-d,1+a-(1+2*a-b-c-d+n),1+a+n},1,k)/
  hyperterm({1+a,1+a-b-c,1+a-b-d,1+a-c-d,1},
  {1+a-b,1+a-c,1+a-d,1+a-b-c-d},1,n),k,n,1);

        (2 a - d + 2 n - c - b) (a + k) (- k + n) (b + k) (c + k) (d + k)
   - -----------------------------------------------------------------------
     n (a + 2 k) (a - b - c - d + n - k) (a - b + n) (a - d + n) (a - c + n)

> sumrecursion(hyperterm({a,1+a/2,b,c,d,1+2*a-b-c-d+n,-n},
  {a/2,1+a-b,1+a-c,1+a-d,1+a-(1+2*a-b-c-d+n),1+a+n},1,k)/
  hyperterm({1+a,1+a-b-c,1+a-b-d,1+a-c-d,1},
  {1+a-b,1+a-c,1+a-d,1+a-b-c-d},1,n),k,n);

                             summ(n) - summ(n - 1)

> hyperrecursion({a,1+a/2,b,c,d,1+2*a-b-c-d+n,-n},
  {a/2,1+a-b,1+a-c,1+a-d,1+a-(1+2*a-b-c-d+n),1+a+n},1,n);

     - (a + n) (a - c - d + n) (a - b - d + n) (a - b - c + n) summ(n - 1)

          + summ(n) (a - d + n) (a - c + n) (a - b + n) (a - b - c - d + n)

> # Gessel-Stanton (6.2)

> WZ(hyperterm({a+1/2,a,b,1-b,-n,(2*a+1)/3+n,a/2+1},
  {1/2,(2*a-b+3)/3,(2*a+b+2)/3,-3*n,2*a+1+3*n,a/2},1,k)/
  hyperterm({(2*a+2)/3,2*a/3+1,(1+b)/3,(2-b)/3,1},
  {(2*a-b)/3+1,(2*a+b+2)/3,2/3,1/3},1,n),k,n);

    6 (a - 1 + 3 n) (a + k) (2 a + 2 k + 1) (- k + n) (b - 1 - k) (b + k)/(

        (a + 2 k) (3 n - k) (3 n - 1 - k) (3 n - 2 - k) (2 a - b + 3 n)

        (2 a + b - 1 + 3 n))

> sumrecursion(hyperterm({a+1/2,a,b,1-b,-n,(2*a+1)/3+n,a/2+1},
  {1/2,(2*a-b+3)/3,(2*a+b+2)/3,-3*n,2*a+1+3*n,a/2},1,k)/
  hyperterm({(2*a+2)/3,2*a/3+1,(1+b)/3,(2-b)/3,1},
  {(2*a-b)/3+1,(2*a+b+2)/3,2/3,1/3},1,n),k,n);

                             summ(n) - summ(n - 1)

> hyperrecursion({a+1/2,a,b,1-b,-n,(2*a+1)/3+n,a/2+1},
  {1/2,(2*a-b+3)/3,(2*a+b+2)/3,-3*n,2*a+1+3*n,a/2},1,n);

     - (3 n - 2 + b) (3 n - 1 - b) (2 a + 3 n) (2 a - 1 + 3 n) summ(n - 1)

          + summ(n) (3 n - 1) (3 n - 2) (2 a - b + 3 n) (2 a + b - 1 + 3 n)

> # The following two sums are identified to be equal, see Strehl (1993)

> sumrecursion(binomial(n,k)^3,k,n);

                 2                   2                                   2
      - 8 (n - 1)  summ(n - 2) - (7 n  - 7 n + 2) summ(n - 1) + summ(n) n

> sumrecursion(binomial(n,k)^2*binomial(2*k,n),k,n);

                 2                   2                                   2
      - 8 (n - 1)  summ(n - 2) - (7 n  - 7 n + 2) summ(n - 1) + summ(n) n

> simplify_combinatorial((binomial(n,k)-binomial(n-2,k))/
  (binomial(n-3,k)-binomial(n-6,k)));

                 (n - 5) (n - 4) (n - 3) (n - 2) (- k + 2 n - 1)
   --------------------------------------------------------------------------
       2                          2
   (3 n  - 24 n - 3 k n + 12 k + k  + 47) (n - 2 - k) (- k + n - 1) (- k + n)

> extended_gosper(binomial(k/2,n),k,2);

                         (1/2 k + 1) binomial(1/2 k, n)
                         ------------------------------
                                      n + 1

> WZ(binomial(n,k)*binomial(k/2,n)*(-1)^k*(-2)^n,k,n,1,2);

                            (- k + n - 1) (- k + n)
                            -----------------------
                            (n - 1) (- k + 2 n - 2)

> extended_sumrecursion(binomial(n,k)*binomial(k/2,n)*(-1)^k*(-2)^n,k,n,1,2);

                             summ(n) - summ(n - 1)

> extended_sumrecursion(binomial(n,k)*binomial(k/2,n)*(-2)^n,k,n,1,2);

                             summ(n) - summ(n - 1)

> hyperrecursion({-n,n+2*a,a},{2/2*a,(2*a+1)/2},2/4,n);

                  (n - 1) summ(n - 2) + (n + 2 a - 1) summ(n)

> hyperrecursion({-n,n+4*a,a},{4/2*a,(4*a+1)/2},4/4,n);

                 - (n - 1) summ(n - 2) + (n + 4 a - 1) summ(n)
\end{verbatim}
}

\label{lastpage}
\end{document}